\documentclass[12pt]{article}

\usepackage{latexsym}
\usepackage{amsmath}
\usepackage{amsfonts}
\usepackage{amssymb}
\usepackage{psfrag}
\usepackage{graphicx}
\usepackage{setspace}

\DeclareGraphicsExtensions{.jpg,.pdf,.eps,.png}

\def\diag{\hbox{diag}}

\def\diag{\hbox{diag}}
\def\log{\hbox{log}}

\def\boxit#1{\vbox{\hrule\hbox{\vrule\kern6pt
          \vbox{\kern6pt#1\kern6pt}\kern6pt\vrule}\hrule}}

\def\bse{\begin{eqnarray*}}
\def\ese{\end{eqnarray*}}
\def\be{\begin{eqnarray}}
\def\ee{\end{eqnarray}}
\def\bq{\begin{equation}}
\def\eq{\end{equation}}
\def\bse{\begin{eqnarray*}}
\def\ese{\end{eqnarray*}}

\newtheorem{Definition}{Definition}

\newtheorem{prop}{Proposition}

\def\part{\partial}

\def\RR{\mathbb{R}}

\newcommand{\bW}{{\bf W}}

\newcommand{\bw}{{\bf w}}

\newcommand{\bU}{{\bf U}}
\newcommand{\bu}{{\bf u}}
\newcommand{\bX}{{\bf X}}
\newcommand{\bx}{{\bf x}}

\newcommand{\bZ}{{\bf Z}}

\newcommand{\bomega}{\mbox{\protect\boldmath $\omega$}}
\newcommand{\btheta}{\mbox{\protect\boldmath $\theta$}}

\newcommand{\balpha}{\mbox{\boldmath $\alpha$}}

\newcommand{\bDelta}{\mbox{\boldmath $\Delta$}}

\newcommand{\bmu}{\mbox{\protect\boldmath $\mu$}}

\newcommand{\bSigma}{\mbox{\protect\boldmath $\Sigma$}}

\newcommand{\bI}{\mbox{\protect\boldmath $I$}}

\newcommand{\by}{{\bf y}}
\newcommand{\bY}{{\bf Y}}

\begin{document}

\renewcommand{\eqref}[1]{(\ref{#1})}
\newcommand{\mb}[1]{\mathbf{#1}}
\newcommand{\mbb}[1]{\mathbb{#1}}
\newcommand{\mt}[1]{\mathrm{#1}}
\newcommand{\rv}{random variable}

\thispagestyle{empty}
\baselineskip=28pt
\vskip 5mm
\begin{center} {\LARGE{\bf Shannon Entropy and Kullback-Leibler Divergence in Multivariate
Log Fundamental Skew-Normal and Related  Distributions }}
\end{center}

\baselineskip=12pt
\vskip 10mm

\begin{center}\large
M. M. de Queiroz \footnote{marinamunizdequeiroz@gmail.com}, R. W. C. Silva \footnote{rogerwcs@ufmg.br} and R. H. Loschi\footnote{loschi@ufmg.br}
\\ 
\small Departamento de Estatística, Universidade Federal de Minas Gerais, Brazil\\

\baselineskip=15pt
\vskip 10mm
\centerline{\today}
\vskip 10mm

\end{center}


\begin{abstract}
This paper mainly  focuses on studying  the Shannon Entropy and Kullback-Leibler divergence
of the multivariate log canonical fundamental skew-normal (LCFUSN) and canonical fundamental skew-normal (CFUSN) families of distributions,
extending previous works.
We relate our results with other well known distributions entropies.
As a byproduct, we also obtain the Mutual Information for distributions in these families. Shannon entropy is used to compare models
fitted to analyze  the USA monthly precipitation data. Kullback-Leibler divergence
is used to cluster regions in Atlantic ocean according to their air humidity level.


\end{abstract}

\section{Introduction}
\setstretch{1}
In recent years,  the interest  in new parametric distributions
able to model  skewness, heavy tails, bimodality and some other data
characteristics that are not well fitted by the usual distributions
is growing. Seeking for flexible and tractable distributions to model 
non-negative data and motivated by some results that recently
appeared in  \cite{SaLoAr13}, the multivariate log-canonical fundamental skew-normal (LCFUSN) family of
distributions is introduced in \cite{QuLoSi15}. Our purpose in this work is to explore some  properties of the LCFUSN family that are useful to solve relevant problems, such as the quantification of the information in a system or
the  definition of optimal designs. These types of problems can be addressed with information theory techniques, such as {\textit {Shannon entropy}} or {\textit{differential entropy}} (see \cite{Sh48})  and Kullback-Leibler (KL) divergence (see \cite{KuLe51, Ku78}). 

The Shannon entropy $H_\bZ$ (hereafter, named entropy)   of a  continuous random vector $\bZ\in\mathbb{R}^n$  can be understood as  the mean  information needed in order to describe the behavior of $\bZ$ whereas the KL divergence measures the inefficiency in assuming that the distribution
is $f_{\bY}$ when the true one is $f_{\bX}$, that is,  it measures the information
lost when $f_{\bY}$ is used to approximate $f_{\bX}$. The KL divergence between $f_{\bX}$  and $f_{\bY}$, 
denoted by $D(f_{\bX}||f_{\bY})$,
is non-negative but it does not define a proper distance measure since it is not symmetric. 
A particular case of the KL divergence is the  {\textit{mutual information}} (MI)  (see \cite{Sh48}) between the random vectors $\bX\in\mathbb{R}^n$ and $\bY\in\mathbb{R}^m$, which is given by  $I_{\bX,\bY}=D(f_{\bX,\bY}(\bx,\by)||f_{\bX}(\bx)f_{\bY}(\by))$.
Therefore, the MI measures the association  between $\bX$ and $\bY$ and is a useful
tool to obtain information about the correlation of such quantities. It also follows  that $I_{\bX,\bY}= I_{\bY, \bX}$
and  $I_{\bX,\bY}=0$  if and only if $\bX$ and $\bY$ are independent. 
The MI can also be obtained  in terms of the  entropies of $\bX$ and $\bY$ through the simple relation 
\begin{equation}\label{propinfo}
I_{\bX,\bY}= H_{\bX}+H_{\bY}-H_{\bX\bY}.
\end{equation}
Additional details on this concepts  can be found in  \cite{CoTh06}.

The entropy of many distributions are already known. In \cite{Ku78} the authors obtained
the entropy of a normal distribution while  in \cite{AhGo89}  the entropies of
several multivariate distributions are derived. In \cite{JaGu08} and \cite{JaGu09} the MI for non normal multivariate location scale families is studied. In
\cite{ArCoGe11} the authors obtained the entropy and MI in
the multivariate elliptical  and skew-elliptical families of distributions
and applied their  results to define an optimal design for an ozone monitoring
station network.  More recently, 
strategies to compute the KL and Jeffreys divergences
in a class of multivariate SN distributions are proposed in \cite{CoAr12}. 

Some information theory concepts have also been considered
in Bayesian Statistics.  For instance, the maximization of the  entropy to
build non informative prior distributions for the parameters is considered in \cite{Ja68}. In
\cite{Be79} the authors 
obtained the so-called reference prior, which corresponds to the Jeffreys
prior in some particular cases, by maximizing the MI 
between the data $\bX$ and the parameter $\btheta$.  In another direction, when the parametric
form of the likelihood can not be identified, 
$I_{\bX,\btheta}$ is minimized in \cite{YuCl99} in order to obtain the minimally
informative likelihood function. In \cite{Mc89} the authors considered the KL divergence to assess the influence of model assumptions in the analysis. A loss function  based on KL divergence to build a criterion for model selection is considered in \cite{Sa02}.
More recently, in \cite{GuMu04} it is  proved that the natural conjugate prior distribution is maximally informative
when it minimizes $I_{\bX,\btheta}$. Also, in \cite{GuSr10} the Bayes estimator for
the entropy by minimizing the expected Bergman divergence is obtained.

In this paper, we mainly  focus on studying  the entropy of the multivariate LCFUSN family of distribution,
as well as  the canonical fundamental skew-normal (CFUSN) distribution defined in \cite{ArGe05}.
As byproducts, we obtain the entropy  of the
multivariate log-skew-normal (LSN) distribution (see \cite{MaGe10}) and generalize
some results obtained in \cite{ArCoGe11} for distributions with
normal kernels.
We obtain the relationship between the entropies of the multivariate LCFUSN and CFUSN distributions
and some related distributions such as
the multivariate normal, multivariate SN (see \cite{AzDa96}) and
the multivariate log-normal (LN) distributions.
We also obtain the KL divergence
to evaluate the ``distance"  between distributions in the CFUSN
family as well as in LCFUSN  family
and also between the LCFUSN  distribution and the LSN distribution.
Some results related to the calculus of the entropy of some univariate
LN and LSN distributions are also shown.
To illustrate our results, we apply
entropy  to compare some log-skewed models fitted to analyze
the USA monthly precipitation data. We also apply the KL divergence
to cluster regions in Atlantic ocean according to their air humidity level.

This work is organized as follows. In Section 2 we present all univariate and multivariate
distributions that are considered in this paper. In Section 3, we find
the entropy of such distributions and obtain the relationship
between them. We
also derive the KL divergence and the MI
for the LCFUSN and CFUSN families of distributions and, as a consequence, for the multivariate
LSN distribution. In Section 4, some issues on Bayesian
inference in the LCFUSN family are discussed and we analyse two procedures
to estimate its entropy and KL divergence. To illustrate our results,
we present two data analysis in Section 5. We also perform an analysis of simulated  data sets
comparing the fitting quality provided by  LCFUSN, LN and LSN  distributions.  In Section 6,
some conclusions and final comments close the paper.

Along the paper, $\phi_n(\bx|\bmu,\bSigma)$ and $\Phi_n(\bx|\bmu,\bSigma)$
denote the pdf  and cumulative distribution function (cdf) of the
multivariate normal distribution $N_n(\bmu,\bSigma)$, respectively.
If $\bmu=0$ such pdf and cdf are denoted by $\phi_n(\bx|\bSigma)$ and $\Phi_n(\bx|\bSigma)$
and, if in addition $\bSigma=\mathbf{I}_n$, we write $\phi_n(\bx)$
and $\Phi_n(\bx)$, respectively. Also, denote by ${\bf{1}}_{n,m}$ and ${\bf{1}}_{n}$
the matrices of ones of order
$n \times m$ and $n\times1$, respectively.

\section{Definitions and Preliminary Results}
The distributions considered in the paper  and
some of its properties are presented in this section.
Although some of them are completely standard, we show them for the benefit of the text.
We begin with the construction of the
LN distribution, since a similar idea is
considered to build  some other log-style  distributions.
It is well-known that if $X \sim N(\mu,\sigma^2)$ and $Y=\exp(X)$ then $Y$ has LN
distribution with parameters $\mu$ and $\sigma^2$, denoted by $Y \sim LN(\mu,\sigma^2)$,
with pdf  given by
$f(y \mid \mu, \sigma)=\frac{1}{y}\phi\left(\ln y|\mu,\sigma\right), \ y \in \mathbb{R}^+.$

Although usual, normality  is not always a reasonable assumption in data analysis if,
for instance, the data have a certain amount of asymmetry or the presence of heavy tails
is realistic. In \cite{Az85} the univariate SN distribution
was introduced by multiplying the pdf of a normal distribution
by a skewing function, which includes an additional parameter to control asymmetry.
We say that $X \sim SN(\mu,\sigma,\alpha)$, where $\mu$, $\sigma^{2}$
and $\alpha \in\mathbb{R}$ are, respectively,  the location, scale and shape parameters,
if its pdf is given by
\begin{equation}
\label{skewnormal}
f(z \mid \mu, \sigma,\alpha)=\frac{2}{\sigma}\phi\left(\frac{z-\mu}{\sigma}\right)\Phi\left(\alpha\left(\frac{z-\mu}{\sigma}\right)
\right), \ z \in \mathbb{R}.
\end{equation}
The normal and the half-normal distributions are particular cases of Equation (\ref{skewnormal}) when $\alpha$ equals zero
and  $\alpha\rightarrow \infty$, respectively.
The expected value and variance of $X \sim SN(\mu,\sigma,\alpha)$ are given, respectively, by
\begin{equation}\label{esperancasn}
E(X)=\mu+\sigma \sqrt{\frac{2}{\pi}} \frac{\alpha}{\sqrt{1+\alpha^2}}\mbox{ and } Var(X)=\sigma^2\left(1-\frac{2\alpha^2}{\pi(1+\alpha^2)}\right).
\end{equation}

Following ideas used in the normal case,
the LSN distribution is introduced in \cite{AzCaKo03}.
Let Z $\sim$ SN($\mu,\sigma,\alpha$) and consider the transformation $Y=\exp(Z)$. Then Y
has the LSN distribution, denoted by $Y\sim LSN(\mu,\sigma,\alpha)$, with
pdf given by
\begin{equation}
\label{logskew}
f(y \mid \mu,
\sigma,\alpha)=\frac{2}{\sigma y}\phi\left(\frac{\ln y-\mu}{\sigma}\right)\Phi\left(\alpha\left(\frac{\ln y-\mu}{\sigma}\right)
\right), y \in \mathbb{R}^+,
\end{equation}
where $\mu \in \mathbb{R}$ is the location parameter, $\sigma^{2}>0$ is the scale parameter and $\alpha \in \mathbb{R}$ is the
shape parameter.
As before, if $\alpha=0$, then Equation (\ref{logskew}) reduces to the LN distribution $LN(\mu,\sigma^2)$.

The multivariate analog of the SN distribution was introduced in  \cite{AzDa96}
and it is  a particular case of the CFUSN family of distributions introduced in \cite{ArGe05}.
We say that  $\bZ$ has a $n$-variate  CFUSN distribution with a $n \times m$ skewness
matrix $\bDelta$, which will be denoted by $\bZ\sim CFUSN_{n,m} (\bDelta)$, if its pdf is given by
\begin{equation}
\label{CFUSN}
f_{\bZ}(\mathbf{z})=2^m\phi_n(\mathbf{z})\Phi_m(\bDelta' \mathbf{z}|{\bf{I}}_m-\bDelta' \bDelta),\,\,\,\,\mathbf{z}\in\mathbb{R}^n,
\end{equation}
where $\bDelta$ is such that ${\bf{I}}_m-\bDelta' \bDelta$ is a positive definite matrix, i.e, $||\bDelta \mathbf{a}|| < 1$, for all unitary vectors $\mathbf{a}\in \mathbb{R}^m$. Here $||^{.}||$
denotes euclidean norm. If  $m=1$ and $\bDelta=(\delta_1,\delta_2,\dots,\delta_n)'$, then we obtain the multivariate SN
family of distributions.
Also, if $m=n$ and $\bDelta=\diag(\delta_1,\delta_2,\dots,\delta_n)$, then Equation
(\ref{CFUSN}) reduces to the product of $n$ SN marginal distributions.
Consequently, for any random sample of the univariate SN distribution  $Y_i\sim SN(\alpha)$, $i=1,\dots,n$,
we have $\bY=(Y_1,\dots,Y_n)' \sim CFUSN_{n,n}(\delta {\bf{I}}_n)$, where
$\delta=\alpha[1+\alpha^2]^{-1/2}$.

In \cite{ArGe05} a location scale version of the CFUSN distribution is introduced. This is accomplished
by considering $\bZ\sim CFUSN_{n,m}(\bDelta)$ and the linear transformation
$\bW=\bmu+\bSigma^{1/2}\bZ$,
where $\bmu$ is the location vector of order $n\times1$ and $\bSigma$ denotes the definite positive scale matrix of
dimension $n \times n$. We say that $\bW\sim CFUSN_{n,m}(\bmu,\bSigma,\bDelta)$ if its pdf is
\allowdisplaybreaks
\begin{eqnarray}
\label{CFUSNLS}
f_{\bW}(\bw)&=&2^{m}|\bSigma|^{-1/2}\phi_{n}(\bSigma^{-1/2}(\bw-\bmu))\nonumber\\
&\times&\Phi_{m}(\bDelta'\bSigma^{-1/2}(\bw-\bmu)|\bI_{m}-\bDelta'\bDelta), \bw \in \mathbb{R}^{n},
\end{eqnarray}
where $|\mathbf{A}|$ stands for $det\,(\mathbf{A})$.

If  data has positive support, the use
of distributions with real support to describe their behavior cannot be appropriate. In the univariate
case there are many different distributions that are useful  to that purpose.
However,  multivariate versions of such univariate distributions are usually intractable.
In \cite{QuLoSi15} the authors introduced the LCFUSN family of distributions in the following way. Let $\bZ=(Z_1, \dots, Z_n)'$ be a $n \times 1$ random
vector and consider the transformations
${\exp({\mathbf{Z}})}=(\exp(Z_1), \dots,\exp(Z_n))'$ and ${\ln {\mathbf{Z}}}=(\ln Z_1, \dots,\ln Z_n)'$.
Hereafter we assume that  $\bDelta^*={\bf{I}}_m-\bDelta'\bDelta$.

\begin{Definition}
\label{Def1}
Let  $\bY$  and $\bZ = \ln \bY$ be $n \times 1$ random vectors. If
$\mathbf{Z} \sim CFUSN_{n,m}(\mathbf{\Delta})$, we say
that $\bY $  has a log-canonical fundamental skew-normal  distribution with skewness matrix 
$\mathbf{\Delta}$ denoted by $\mathbf{Y} \sim LCFUSN_{n,m}(\mathbf{\Delta})$.
The pdf  of $\mathbf{Y}$ is
\begin{equation}\label{LCFUSN}
f_{\mathbf{Y}}(\mathbf{y})=2^{m}\left(\prod_{i=1}^{n}y_{i}\right)^{-1}
\phi_{n}(\ln \mathbf{y})\Phi_{m}(\mathbf{\Delta}'\ln \mathbf{y}|\bDelta^*),
\,\,\,\,\mathbf{y} \in \mathbb{R}^{n^{+}},
\end{equation}
where $\mathbf{\Delta}$ is a $n \times m$ matrix such that $||\mathbf{\Delta} \textbf{a}||<1$, for all unity vectors
\textbf{a} $\in \mathbb{R}^{m}$.
\end{Definition}

The location-scale version of the LCFUSN distribution is obtained
by considering a $n\times 1$ random vector $\bW \sim CFUSN_{n,m}(\bmu,\bSigma, \bDelta)$ and the
transformation $\mathbf{U}=\exp(\bW)$. In this case, the pdf of  $\mathbf{U}\sim LCFUSN_{n,m}(\bmu,\bSigma, \bDelta)$ is
\allowdisplaybreaks
\begin{equation}
\label{fdpescalaloc}
f_{\bU}(\bu)=\frac{2^{m}|\bSigma|^{-1/2}}{\prod_{j=1}^{n}u_{j}}\phi_{n}(\bSigma^{-1/2}(\ln \bu-\bmu))
\Phi_{m}(\bDelta'\bSigma^{-1/2}(\ln \bu-\bmu)|\bDelta^*),
\end{equation} for all $\mathbf{u} \in \mathbb{R}^{n^{+}}$, where $\bmu$
is a $n\times1$ location vector, $\bSigma$ a $n\times n$ definite positive scale matrix and
$\bDelta$ a skewness matrix.
 The LCFUSN family of distributions generalizes the multivariate
LSN family defined in \cite{MaGe10}, which can be obtained from Equation (\ref{LCFUSN}) by taking $m=1$ and $\balpha'=(1-\mathbf{\Delta}'
\mathbf{\Delta})^{-\frac{1}{2}}\mathbf{\Delta}'{\bSigma}^{-\frac{1}{2}}{\bomega}$, where ${\bomega}=diag({\bSigma})^{\frac{1}{2}}$ is a $n\times n$ scale matrix. 

%


\section{Information Theory in the CFUSN  and LCFUSN Distributions }
\label{Se3}

In this section we  calculate  the entropy and
the KL divergence  for the multivariate
CFUSN  and the multivariate LCFUSN families of distributions and associate them to
the entropy of related distributions. We start with the   univariate case and obtain
the entropy of the LSN distribution (see \cite{AzCaKo03})   and
its relationship with the normal, SN and LN entropies.

\subsection{Univariate Cases}

It is  well-known that the maximum entropy  among  all univariate continuous symmetric
distributions  is observed for  the normal distribution (see \cite{CoTh06}).
If  $X \sim N(\mu,\sigma^2)$, then the entropy of $X$ is
\begin{equation}
\label{entnormal2}
H_{N(\mu,\sigma^2)}=\frac{1}{2}\ln \sigma^{2} + \frac{1}{2}(1+\ln(2\pi)),
\end{equation}
which increases with the variance of the distribution.
Opposed to what is observed for $H_{N(\mu,\sigma^2)}$,
the entropy of a LN distribution also depends on the location parameter $\mu$. Formally, if
$X\sim LN(\mu,\sigma^{2})$, then
\begin{equation}\label{entrolog}
H_{LN(\mu,\sigma^{2})}= H_{N(\mu,\sigma^{2})} + \mu.
\end{equation}
The entropy of a LN distribution is an increasing function of $\mu$ and is higher
than the entropy of  $N(\mu,\sigma^2)$ if and only if $\mu>0$.
In \cite{ArCoGe11} the entropy of the multivariate
skew-elliptical class of distributions is obtained.
Particularly, they show that  the entropy of $X\sim SN(\mu,\sigma^{2},\alpha)$, 
a special case of such a class, is a function of the entropy of the $N(\mu,\sigma^2)$
distribution and is given by
\begin{equation}
\label{entroskew}
H_{SN(\mu,\sigma^{2}, \alpha)}= H_{N(\mu,\sigma^{2})}-E_{X_0}[\ln (2\Phi(\alpha X_{0}))],
\end{equation}
where $X_{0}\sim SN(\alpha)$.
As expected, if $\alpha=0$, we have that $H_{SN(\mu,\sigma^2,\alpha)}=H_{N(\mu,\sigma^2)}$.
Moreover, it follows from  the monotone convergence theorem that
\begin{equation}
\label{limitesn}
\lim_{\alpha\to \pm  \infty}H_{SN(\mu,\sigma^2,\alpha)}=H_{N(\mu,\sigma^2)}-\ln(2).
\end{equation}
%

As far as we know, the next result is new and provides the entropy of the LSN distribution.
\begin{prop}
\label{entropialogskewuni}
If $Y\sim LSN(\mu,\sigma^{2},\alpha)$  then the entropy of $Y$ is
\begin{equation}\label{entrologskew}
H_{LSN(\mu,\sigma^{2},\alpha)}= H_{SN(\mu,\sigma^{2},\alpha)}+E(X),
\end{equation}
where $X \sim SN(\mu,\sigma^{2},\alpha)$ and $E(X)$ is given in Equation (\ref{esperancasn}).
\end{prop}

The proof of Proposition \ref{entropialogskewuni} follows by noticing that $\ln f_Y(Y)\buildrel d \over =-X+\ln f_X(X)$, where $X \sim SN(\mu,\sigma^{2},\alpha)$.
Besides, the relationship between the entropies of the LSN, LN and normal distributions
follows immediately from results in Equations (\ref{entrolog}) and  (\ref{entroskew}).
The entropy in Equation $(\ref{entrologskew})$ can be written as 
$ H_{LSN(\mu,\sigma^{2}, \alpha)} = H_{N(\mu,\sigma^{2})}-E[\ln (2\Phi(\alpha X_{0}))] + E(X)$
and  
$  H_{LSN(\mu,\sigma^{2}, \alpha)} = H_{LN(\mu,\sigma^{2})}-E[\ln (2\Phi(\alpha X_{0}))]+\mu+E(X),$
where   $X \sim SN(\mu,\sigma^{2},\alpha)$ and  $X_{0}\sim SN(\alpha)$. Moreover, from Equations  (\ref{esperancasn}) and (\ref{limitesn}), 
\allowdisplaybreaks
\begin{eqnarray*}
\lim_{\alpha\to{\pm \infty}}H_{LSN(\mu,\sigma^2,\alpha)} &=& H_{N(\mu,\sigma^2)}+\mu \pm \sigma\sqrt{\frac{2}{\pi}}-\ln(2).
\end{eqnarray*}
Also, if $\alpha =0$ in Equation (\ref{entrologskew}), then 
$H_{LSN(\mu,\sigma^2,\alpha)}=H_{LN(\mu,\sigma^2)}.$
\begin{figure}[htb]
  \centering
  \includegraphics[width=4.3cm]{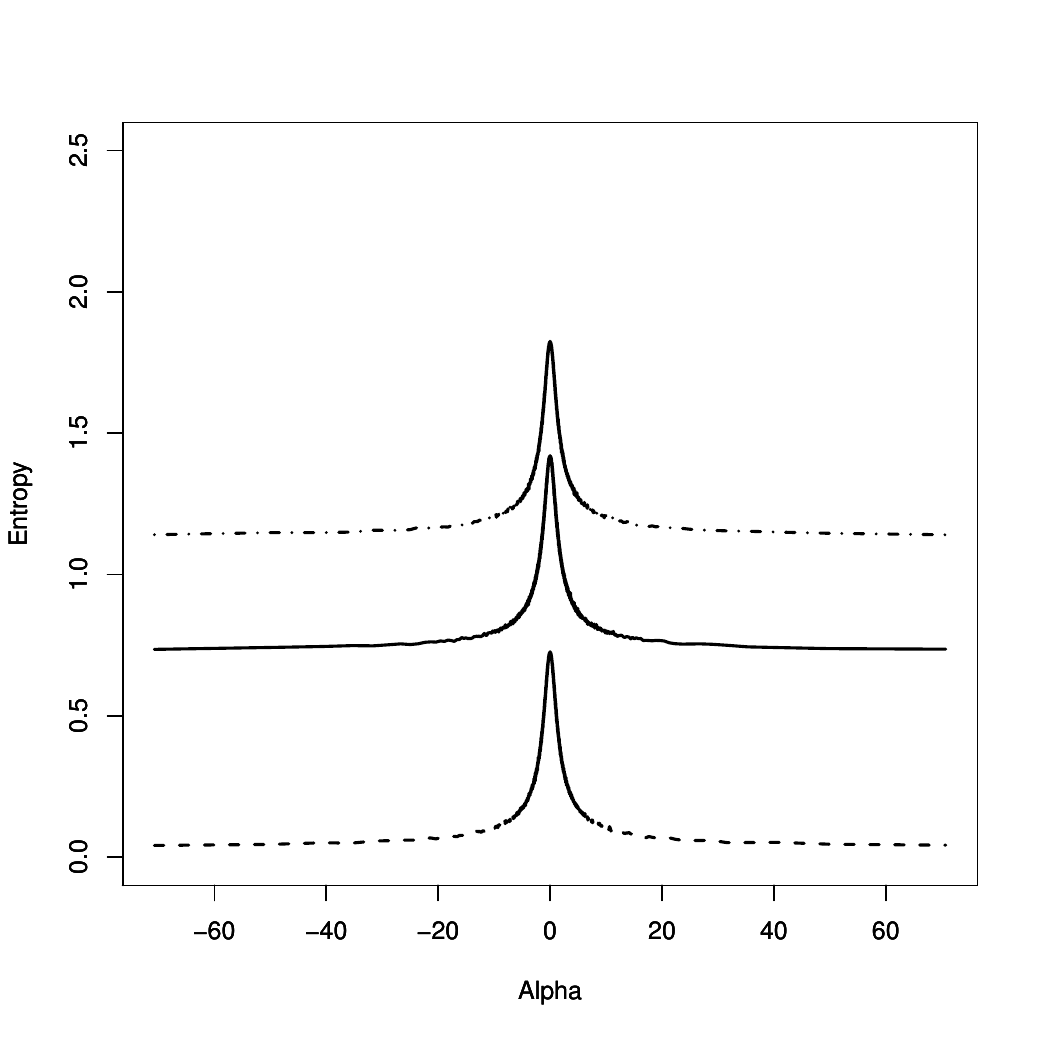}
  \includegraphics[width=4.3cm]{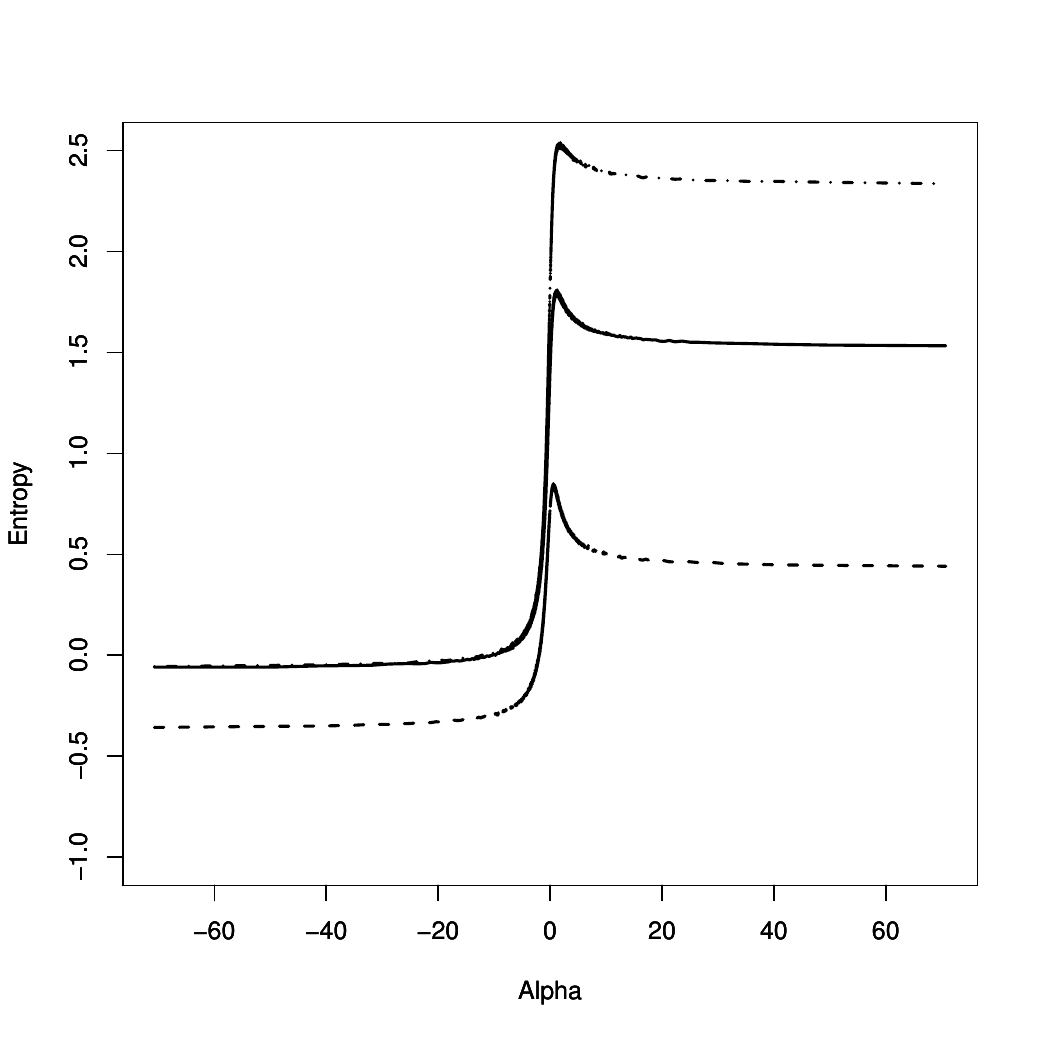}\\
  \caption{Entropy of the $SN(0,\sigma^2,\alpha)$ (left) and $LSN(0,\sigma^2,\alpha)$ (right)
  for $\sigma^2 =0.5$ (dotted line), $1$ (solid line) and $1.5$(dotted-dashed line).}
  \label{entropia}
\end{figure}

Figure \ref{entropia} displays the entropies of $SN(0,\sigma^2,\alpha)$  and $LSN(0,\sigma^2,\alpha)$
as a function of the skewness parameter $\alpha$. The expectations involved in their
computation were approximated using Monte Carlo methods.
Figure \ref{entropia} suggests that $H_{LSN(\mu,\sigma^2,\alpha)}$ is an increasing function of $\alpha$ until a
global maximum $H_{LSN(\mu,\sigma^2,\alpha^*)}$ and decreasing afterwards. A similar behavior is observable for $H_{SN(\mu,\sigma^2,\alpha)}$. Also, we note that the curve of $H_{SN(\mu,\sigma^2,\alpha)}$ has a symmetric shape and if we increase $\sigma^2$ by one unity it is translated by a factor
of $\ln(\sigma^2 +1)$.

\subsection{Multivariate cases}

In this section we obtain the entropies of the CFUSN and LCFUSN distributions
and their connections with the entropies of a multivariate normal, multivariate SN
and multivariate LSN (see \cite{AzDa96} and \cite{MaGe10}) distributions. We
extend some results obtained in \cite{ArCoGe11} for distributions
with normal kernel.

Let $\bZ$ be a $n$-dimensional random vector
such that $\bZ \sim N_{n}(\bmu, \bSigma)$. The entropy of   $\bZ $ is
\begin{equation}\label{entropianorm}
H_{N_{n}(\bmu, \bf{\Sigma})}=\frac{1}{2}\ln |\bSigma| + \frac{n}{2}(1+\ln(2\pi))=\frac{1}{2}\ln |\bSigma|  +
H_{N_{n}({\bf{0}}, {\bf{I}}_{n})},
\end{equation}
where $H_{N_{n}({\textbf{0}}, {\bf{I}}_{n})}= \frac{n}{2}(1+ \ln(2\pi)) $ is the entropy of a random vector with
standard $n$-variate normal distribution.
As in the univariate case, $H_{N_{n}(\bmu, \bf{\Sigma})}$ does not depend on the location parameter $\bmu$.

The CFUSN family 
is a more general class of skewed distributions with normal kernel. In \cite{ArGe05} many of its properties are obtained.
In particular, it is shown  that if $\bX\sim CFUSN_{n,m}(\bmu,\bSigma,\mathbf{\Delta})$ then
\begin{equation}
\label{EspCFSUN}
E(\bX)=\mathbf{\bmu}+\sqrt{\frac{2}{\pi}}\bSigma^{1/2}\bDelta\mathbf{1}_m,\,\,\,\,\,\,Var(\bX)=\bSigma-\frac{2}{\pi}\bSigma^{1/2}\bDelta\bDelta'\bSigma^{1/2}.
\end{equation}
They also prove that  if $(\bX_1,\dots,\bX_n)$ and $(\bDelta_1, \dots, \bDelta_n)$ are partitions of $\bX\sim CFUSN_{n,m}{(\bDelta)}$
and $\bDelta$ respectively, then 
$X_i \sim CFUSN_{1,m}(\bDelta_i)$, $i=1,\dots,n$.
Despite their great contribution, the authors in  \cite{ArGe05} do not obtain results related to the
entropy in the CFUSN family. In the next proposition we provide an expression for this entropy. Proof is given in Appendix.
\begin{prop}\label{enlsncap2}
If  $\bX \sim CFUSN_{n,m}(\bmu,\bSigma, \bDelta)$, then the entropy of the canonical fundamental skew normal random vector
$\bX$ is
\allowdisplaybreaks
\begin{eqnarray}
\label{entropicfsunm}
H_{CFUSN(\bmu, \bf{\Sigma}, \bf{\Delta})}&=&H_{N_{n}(\bf{0}, \bf{I_{n}})}+\frac{1}{2}\ln |\bSigma| +\frac{1}{\pi}[\mathbf{1}_m^{'}\bDelta^{'}\bDelta\mathbf{1}_m-tr(\bDelta\bDelta^{'})]\\\nonumber
&-&E_{\bX_{0}}[\ln(2^{m}\Phi_{m}(\bDelta'\bX_{0}|\bDelta^{*}))]
\end{eqnarray}
where $\bX_{0} \sim CFUSN_{n,m}(\bDelta)$.
\end{prop}
Some interesting results are obtained for particular structures of the
scale $\bSigma$ and the skewing  $\bDelta$ matrices. If $\bSigma$ is a covariance matrix, the relationship between
$H_{N_{n}(\bmu, \bf{\Sigma})}$ and $H_{CFUSN(\bmu, \bf{\Sigma}, \bf{\Delta})}$
follows from Equations (\ref{entropianorm}) and (\ref{entropicfsunm}) . If in addition 
$\bDelta\bDelta^{'}$ is  diagonal, then the entropy in Equation (\ref{entropicfsunm}) is simplified to 
\allowdisplaybreaks
\begin{eqnarray}\label{entropcfusn}
  H_{CFUSN_{n,m}(\bmu, \mathbf{\Sigma}, \mathbf{\Delta})}=H_{N_{n}(\bmu, \bf{\Sigma})}
  -E[\ln(2^{m}\Phi_{m}(\bDelta'\bX_{0}|\bDelta^*))],
\end{eqnarray}
where $\mathbf{X}_{0} \sim CFUSN_{n,m}(\bDelta)$.
The result in Equation (\ref{entropcfusn}) generalizes those obtained in \cite{ArCoGe11} for distributions
with normal kernel, including  the multivariate SN distribution defined in
\cite{AzDa96}. 

\begin{figure}[htb]
  \centering
   \includegraphics[width=4.3cm]{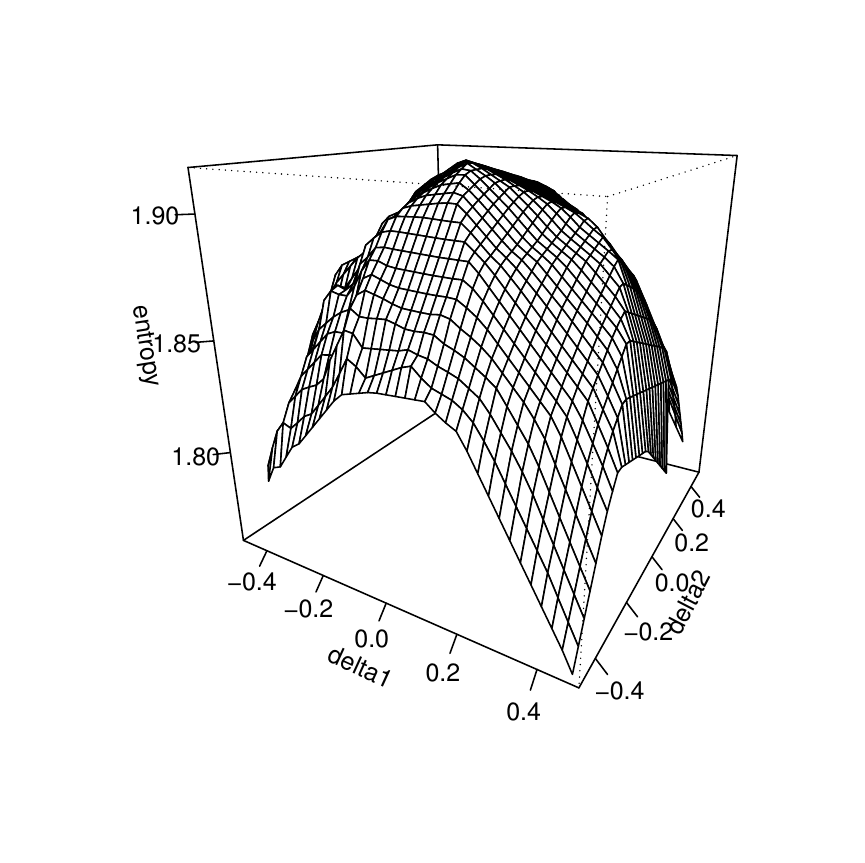}
   \includegraphics[width=4.3cm]{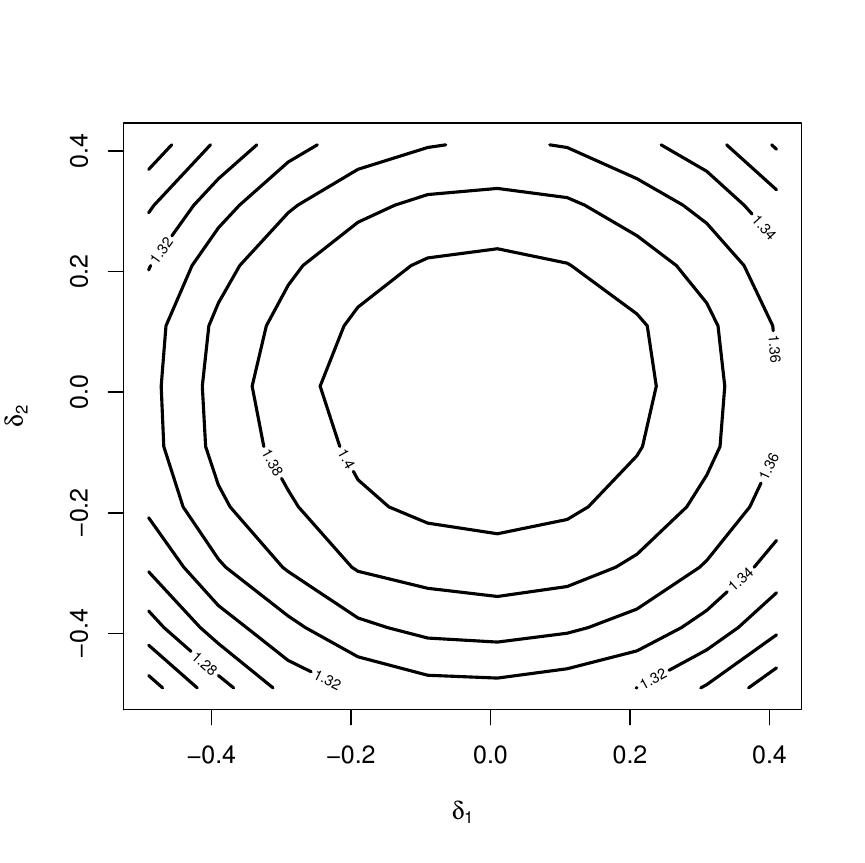}
   \includegraphics[width=4.3cm]{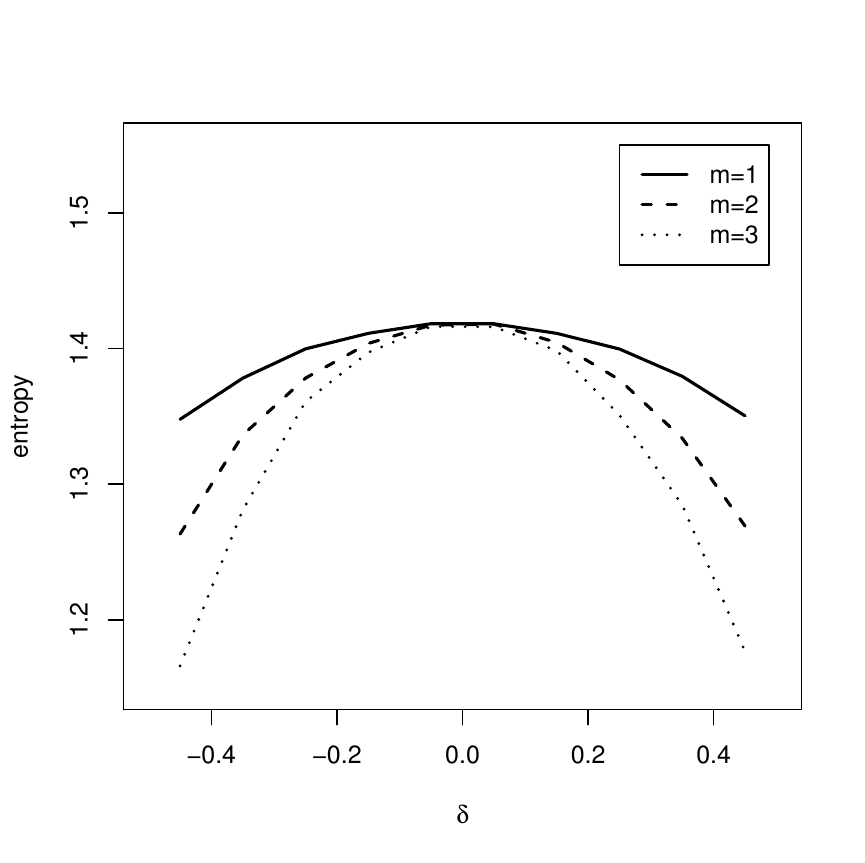}
  \caption{Entropy of the $CFUSN_{1,2}(\bDelta)$ with $\bDelta=(\delta_1,\delta_2)$ (left and middle) and
  the $CFUSN_{1,2}(\delta {\bf{1}}_{m})$(right).}
  \label{EntCFUSN}
\end{figure}

Figure \ref{EntCFUSN} shows the entropy of the standard CFUSN family in the multivariate and univariate cases. We observe
that the entropy is concave and presents symmetric behavior around zero.
The maximum entropy is obtained when the skewing parameter
is zero, that is, in the normal case. Besides, 
the entropy decay is smooth  for small values of $m$. We also note that, for fixed
$\bDelta$, the smaller the value of $m$, the higher the entropy. Finally, the entropy tends
to be closer to the normal entropy for
values of $\bDelta$ around zero.

In Proposition \ref{propentropialsnfc} we obtain  the entropy of the LCFUSN distribution introduced in \cite{QuLoSi15} which pdf is given in Equation (\ref{LCFUSN}).

\begin{prop}
\label{propentropialsnfc}
If $\bZ \sim LCFUSN_{n,m}(\bmu, \bSigma, \bDelta)$, then the entropy of $\bZ$  is
\begin{equation}
\label{entropia2}
H_{LCFUSN_{n,m}(\bmu, \bf{\Sigma}, \bf{\Delta})}=H_{CFUSN_{n,m}(\bmu, \bf{\Sigma}, \bf{\Delta})}+ \sum_{i=1}^{n} E(X_{i}),
\end{equation}
where $X_{i}$ is the $ith$ component of the vector $\bX \sim CFUSN_{n,m}(\bmu, \bSigma, \bDelta)$.
\end{prop}

The proof of Proposition \ref{propentropialsnfc} is straightforward by noticing  
that  if  $\bX \sim CFUSN(\bmu,\bf{\Sigma},\bf{\Delta})$ then
$\ln f_{\bY}(\bY)\buildrel d \over =\ln f_{\bX}(\bX) -{\bf{1}}_n^{'} \bX$.
We can also rewrite the entropy  in Equation (\ref{entropia2}) considering the
results  in  Equation (\ref{EspCFSUN}) and Proposition $\ref{enlsncap2}$, obtaining 
\allowdisplaybreaks
\begin{eqnarray}\label{entropiiiiiaaa}
&&H_{LCFUSN_{n,m}(\bmu, \bf{\Sigma}, \bf{\Delta})}= H_{N_{n}(\bf{0}, \bf{I_{n}})} +\frac{1}{2} \ln |\bSigma|+\frac{1}{\pi}[\mathbf{1}_m^{'}\bDelta^{'}\bDelta\mathbf{1}_m -tr(\bDelta\bDelta^{'})]\nonumber\\
 &+&\sum_{i=1}^{n} \mu_{i}+\sqrt{\frac{2}{\pi}}{\bf 1}_n' \bSigma^{1/2}\bDelta\textbf{1}_{m} -E_{X_{0}}[\ln(2^{m}\Phi_{m}(\bDelta'\bX_{0}|\bDelta^{*}))],
\end{eqnarray}
where the random vector $\bX_{0} \sim CFUSN_{n,m}(\bDelta)$.
If $\bSigma$ is a covariance matrix, the relationship between the LCFUSN and the non standard multivariate normal entropies follows from Equations  (\ref{entropianorm}) and (\ref{entropiiiiiaaa}) . If additionally $\bDelta^{'}\bDelta$ is a diagonal matrix then we obtain
\allowdisplaybreaks
\begin{eqnarray}
\label{entroplcfusn}
  H_{LCFUSN_{n,m}(\bmu, \bf{\Sigma}, \mathbf{\Delta})}&=&H_{N_{n}(\bmu, \mathbf{\Sigma})}+\sum_{i=1}^{n} \mu_{i}
+ \sqrt{\frac{2}{\pi}}{\bf 1}_n' \bSigma^{1/2}\bDelta\textbf{1}_{m} \nonumber\\&-&E_{\bX_0}[\ln(2^{m}\Phi_{m}(\bDelta'\bX_{0}|\bDelta^*))],
\end{eqnarray}
where $\mathbf{}X_{0} \sim CFUSN_{n,m}(\bDelta)$ and $\mu_{i}$ is the $i$th component of  vector $\bmu$.



The entropies of the  multivariate LSN and the multivariate LN distributions are obtained easily from Proposition \ref{propentropialsnfc}, since these distributions are special cases of the LCFUSN distribution.

\begin{figure}[h!]
  \centering
   \includegraphics[width=4.3cm]{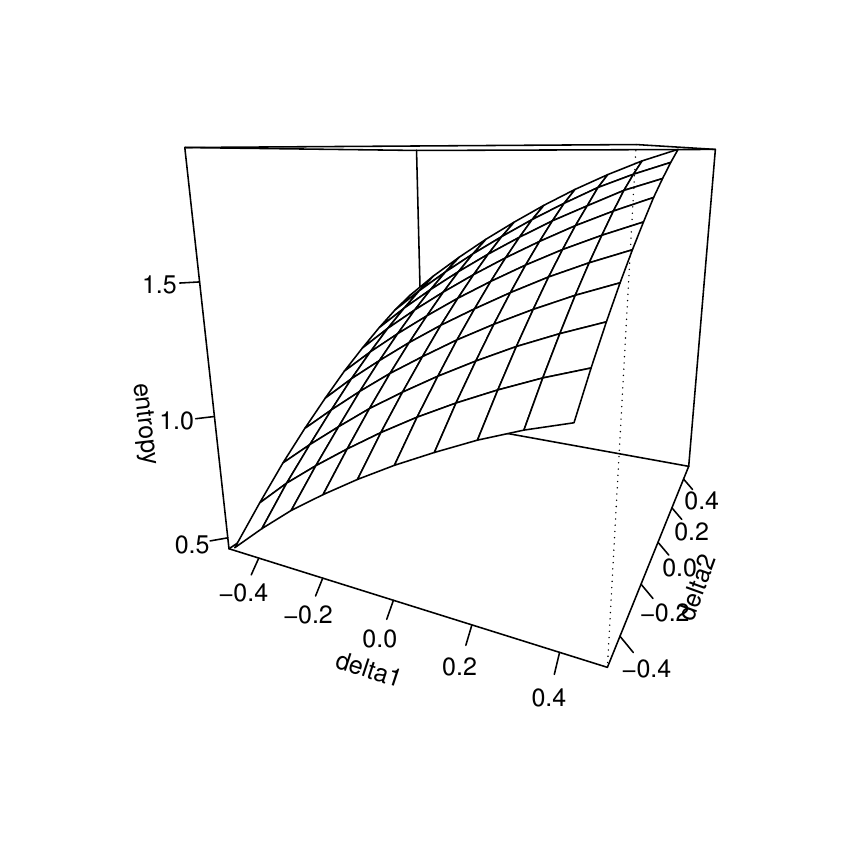}
    \includegraphics[width=4.3cm]{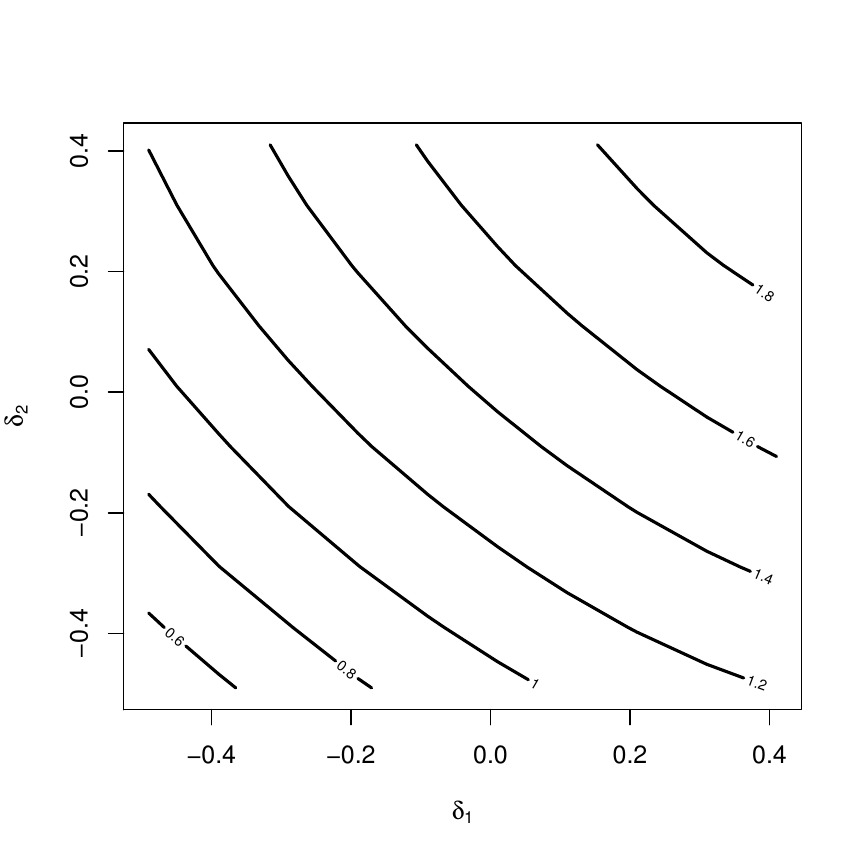}
    \includegraphics[width=4.3cm]{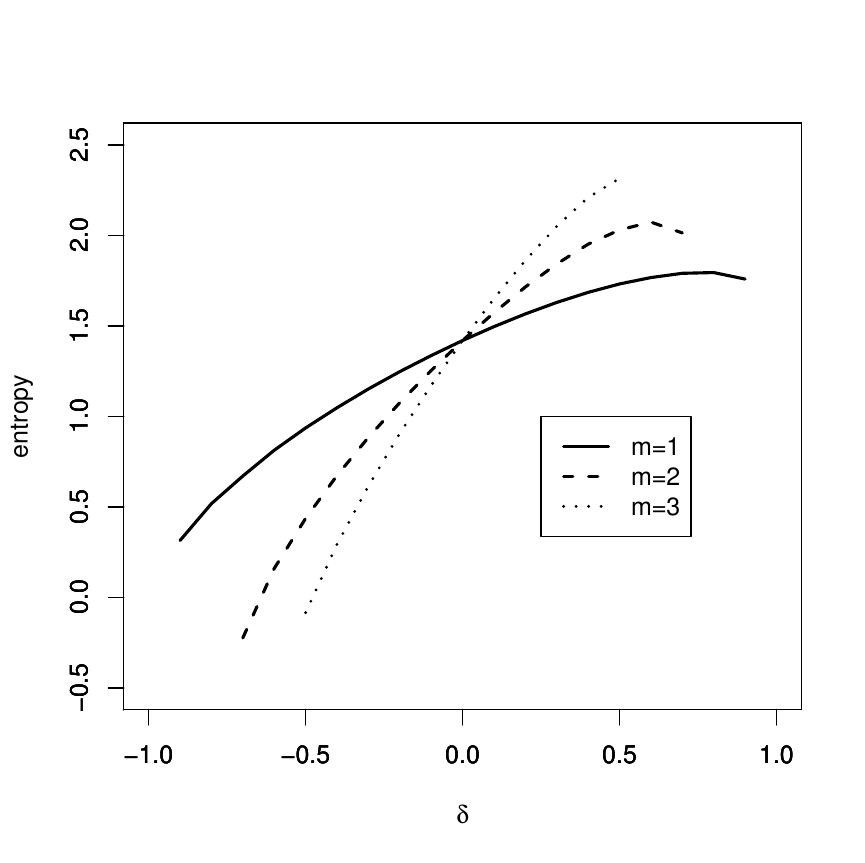}
   \caption{Entropy of the $LCFUSN_{1,2}(\bDelta)$ with $\bDelta=(\delta_1,\delta_2)$ (left and middle) and
   $\bDelta=$ $\delta{\bf{1}}_{m}$ (right).}
  \label{MEnLcfusn}
\end{figure}

Figure  \ref{MEnLcfusn}  shows  the entropy of the univariate  standard  $LCFUSN_{1,2}(\bDelta)$ 
as a function of  $\bDelta$. We assume two different structures for the skewness matrix $\bDelta$.
The plot on the right  also shows the influence of $m$ in the entropy of a
$LCFUSN_{1,m}(\delta{\bf{1}}_{m})$.  It can be seen that the entropy of the LCFUSN
increases  with $\delta$ and is smooth for small values of $m$. Moreover, for $\delta <0$, the highest the value of $m$, the smallest the entropy. The opposite is
observed if $\delta>0$.

\subsection{KL Divergence and MI in the CFUSN and LCFUSN Families of Distributions}
\label{MIKLD}
An useful tool  to compare two distributions is the so called KL divergence.
This quantity measures the inefficiency of assuming that the true distribution is $f_{\bX}$ whereas it is $f_{\bY}$. In the next proposition we obtain the KL divergence in the CFUSN and LCFUSN families. It is remarkable that Proposition \ref{propentropiarelativa} provides the KL divergence for all families of distributions considered in this work. As can be seen below, the KL divergence is invariant under the exponential transformation.
\begin{prop}
\label{propentropiarelativa} In the following cases:
\begin{itemize}
\item[(i)]if $\bZ \sim CFUSN_{n,m_1}(\bmu_1, \bSigma_1, \bDelta_1)$ and $\bY \sim CFUSN_{n,m_2}(\bmu_2, \bSigma_2, \bDelta_2)$ and
\item[(ii)] if $\bZ \sim LCFUSN_{n,m_1}(\bmu_1, \bSigma_1, \bDelta_1)$ and $\bY \sim LCFUSN_{n,m_2}(\bmu_2, \bSigma_2, \bDelta_2)$,
\end{itemize}
the KL divergence between $\bY$ and $\bZ$ is given by
\allowdisplaybreaks
\begin{eqnarray}
\label{entrorela}
D(f_{\bY}||f_{\bZ})&=& \frac{1}{2} \ln\left[\frac{|\bSigma_{1}|}{|\bSigma_{2}|} \right] -\frac{n}{2}
+E_{\bY}\left( \ln \left[\frac{2^{m_{2}-m_{1}}\Phi_{m_{2}}(\bDelta_{2}'\bSigma_{2}^{-\frac{1}{2}}( \bY -\bmu_{2})|\bDelta_{2}^{*})}
{\Phi_{m_{1}}(\bDelta_{1}'\bSigma_{1}^{-\frac{1}{2}}( \bY -\bmu_{1})|\bDelta_{1}^{*})} \right]\right)\nonumber \\
&+& \frac{1}{2}[(\bmu_1-\bmu_2)'\bSigma_1^{-1}(\bmu_1-\bmu_2)+tr(\bSigma_1^{-1}\bSigma_2)]\nonumber\\
&+&\frac{1}{\pi}[{\bf{W}}'{\bSigma}_1{\bf{W}}-tr(\bSigma_2^{-1/2}{\bDelta}_2{\bDelta}_2^{'}{\bSigma}_2^{-1/2})-\mathbf{1}_{m_2}^{'}\bDelta_2^{'}\bDelta_2\mathbf{1}_{m_2}+tr(\bDelta_2\bDelta_2^{'})]\nonumber\\
&+&\frac{1}{\sqrt{2\pi}}[({\bmu}_2'-{\bmu}_1'){\bf{W}}+{\bf{W}}'({\bmu}_2-{\bmu}_1)],
\end{eqnarray}
where ${\bf{W}}=\bSigma_1^{-1}\bSigma_2^{-1/2}\bDelta_2\mathbf{1}_{m_2}$ and $\bDelta_{i}^{*}= {\bf{I}}_{m_i} -\bDelta_i' \bDelta_i$.
\end{prop}

The proof of item $(i)$ in Proposition \ref{propentropiarelativa} can be found in the appendix. Item $(ii)$ follows  by observing that if $\bZ \sim CFUSN_{n,m_1}(\bmu_1, \bSigma_1, \bDelta_1)$, $\bY \sim CFUSN_{n,m_2}(\bmu_2, \bSigma_2, \bDelta_2)$, $\bf U=\exp{(\bZ)}$ and $\bf V=\exp{(\bY)}$, then 
\allowdisplaybreaks
\begin{eqnarray}
D(f_{\bU}||f_{\bf V})&=&E\left(\log\frac{f_{\bf U}(\bU)}{f_{\bf V}(\bU)}\right)=E\left(\log\frac{f_{\bf Z}(\log\, \bU)}{f_{\bf Y}(\log\, \bU)}\right)\nonumber\\
&=&E\left(\log\frac{f_{\bf Z}(\bf Z)}{f_{\bf Y}(\bf Z)}\right)=D(f_{\bf Z}||f_{\bf Y }).\nonumber
\end{eqnarray}
%

The KL divergence between two $n$-variate normal distributions is obtained from
Equation (\ref{entrorela}) by assuming $\bDelta_1$ and $\bDelta_2$ as null matrices.
Proposition \ref{propentropiarelativa} also provides the KL divergence between the  LSN 
and  LCFUSN distributions. If we consider the parametrization of the LSN distribution assumed in \cite{MaGe10}, then Equation (\ref{entrorela}) becomes
\begin{equation*}
D(f_{\bY}||f_{\bZ})=E_{\bX_{0}}\left[ \ln \left(\frac{ 2^{m-1}\Phi_{m}(\bDelta\bX_{0}|\bDelta^{*})}{ \Phi(\balpha'\bomega^{-1}\bSigma^{\frac{1}{2}}\bX_{0})}\right)\right],
\end{equation*}
where $\bX_{0} \sim CFUSN_{n,m}(\bDelta)$, $\balpha\in\mathbb{R}^n$ and ${\bomega}=diag(\bSigma)^{1/2}$.

\begin{figure}[h!]
  \centering
   \includegraphics[width=4.5cm]{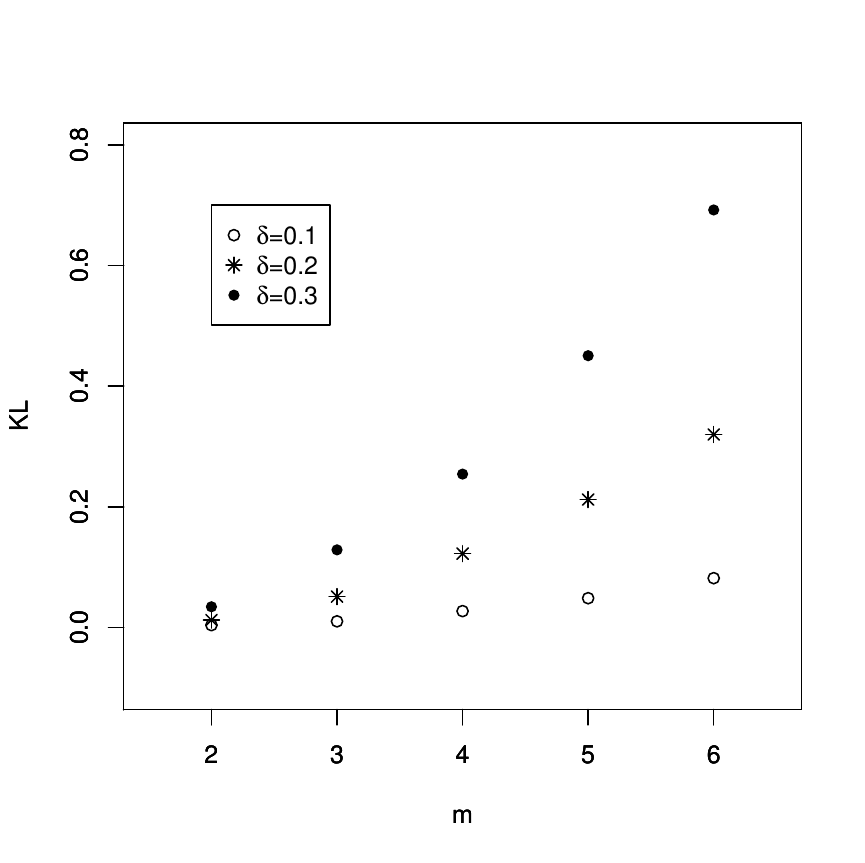}
    \includegraphics[width=4.5cm]{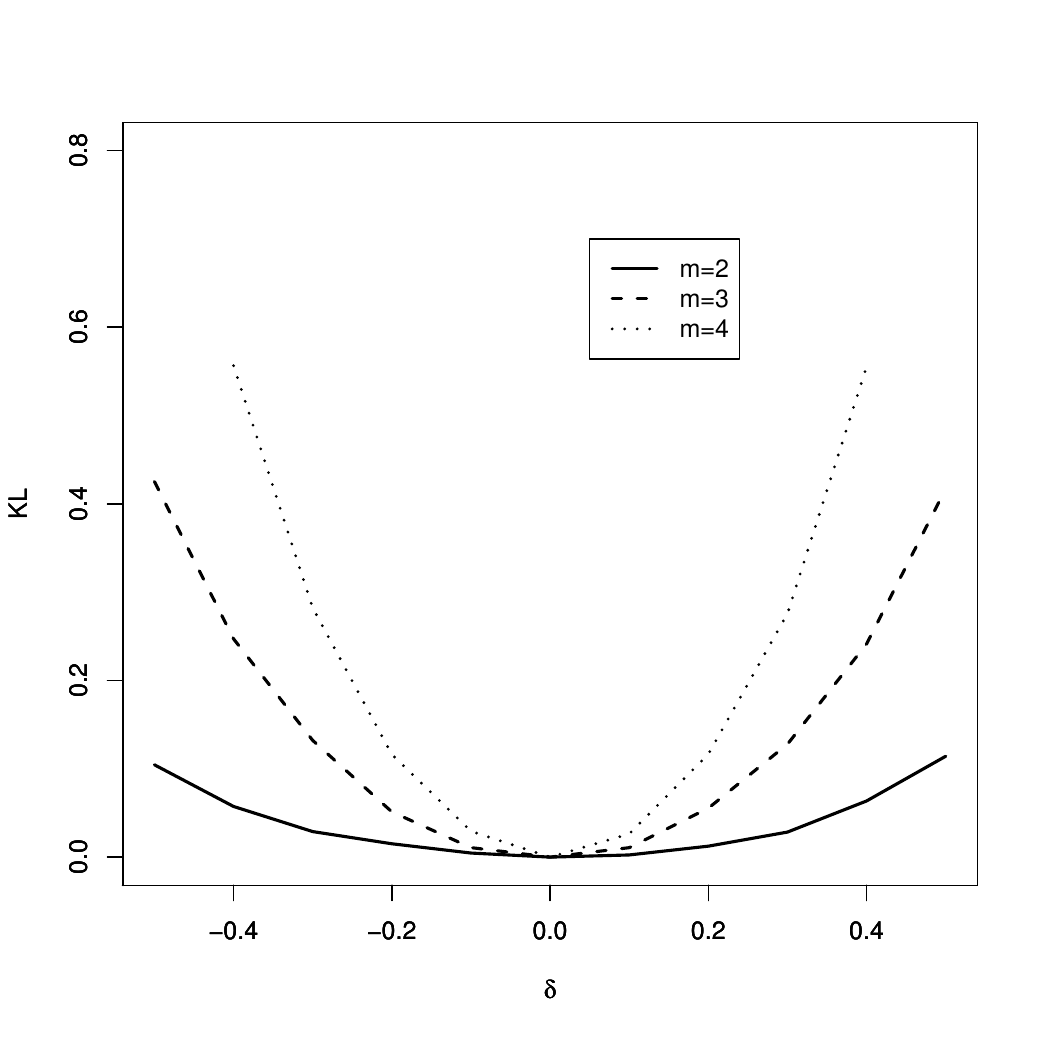}
   \caption{KL divergence between $LCFUSN_{1,m}(\bDelta)$ and $LSN_1(0,1,\alpha)$ with $\balpha'=(1-\bDelta'\bDelta)^{-\frac{1}{2}}\bDelta'$  and
   $\bDelta=$ $\delta{\bf{1}}_{m}$ .}
  \label{KL}
\end{figure}

Graphics displayed in Figure \ref{KL} disclose that, for fixed $\delta$,  the KL divergence between  $LCFUSN_{1,m}(\bDelta)$ and $LSN_1(0,1,\alpha)$
increases as $m$ increases. Similar behavior is observed  for fixed $m$ when $\delta$ is increasing and positive. Moreover, the KL divergence seems
to be symmetric in $\delta$ at least in cases  $m=2,3 \mbox{ and }4$.

Let $\bX_1$ and $\bX_2$ be a partition of $\bX$.
Suppose we want to quantify the amount of information $\bX_1$ brings about
$\bX_2$ when $\bX$ has distribution in the CFUSN or the  LCFUSN family. Since these families are closed under marginalization (see \cite{QuLoSi15, ArGe05}), from  Equation (\ref{propinfo}) we obtain 
\begin{equation}\label{infolcfusn}
I_{\bX_{1}\bX_{2}}= E_{X}\left[ \ln \left(\frac{\Phi_{m}(\bDelta'\bX|\mathbf{I}_{m}-\bDelta'\bDelta)}{2^{m}
\Phi_{m}(\bDelta_{1}'\bX_{1}|\mathbf{I}_{m}-\bDelta_{1}'\bDelta_{1})
\Phi_{m}(\bDelta_{2}'\bX_{2}|\mathbf{I}_{m}-\bDelta_{2}'\bDelta_{2})} \right)\right],
\end{equation}
in the following cases:
\begin{itemize}
\item[(i)] if $\bX \sim CFUSN_{n,m}(\bDelta)$, 
$\bX_{1} \sim CFUSN_{n_{1},m}(\bDelta_{1})$ and $\bX_{2} \sim CFUSN_{n_{2},m}(\bDelta_{2})$; 
\item[(ii)]if $\bX \sim LCFUSN_{n,m}(\bDelta)$,  
$\bX_{1} \sim LCFUSN_{n_{1},m}(\bDelta_{1})$ and $\bX_{2} \sim LCFUSN_{n_{2},m}(\bDelta_{2})$.
\end{itemize}

Observe that the MI on these families of distributions is also invariant under the exponential transformation, since it is a particular case of the KL divergence.


\section{Bayesian Estimation of the LCFUSN Entropy}
\label{SeBaE}

The  entropy and the KL divergence depend on the parameters of the LCFUSN distribution and thus, under the
Bayesian paradigm, are random  quantities. In order to estimate  such quantities it is necessary
to obtain their posterior distributions, which can be a hard task
if we search for closed expressions for them. However, good approximations
can be obtained using MCMC methods.  We only take into consideration
particular univariate cases of the LCFUSN family since one of our goals in Section \ref{Apl}
is to evaluate the effect of increasing  $m$ in data fitting.
Similar strategy can be used in the multivariate case.

To achieve our goal,  let $Y_{1},..., Y_{L}$ be a random sample of $Y \mid \mu, \sigma, \bDelta \buildrel iid
\over\sim LCFUSN_{1,m}(\mu,\sigma^2,\bDelta)$,
which induces the following likelihood function:
\allowdisplaybreaks
\begin{equation}
\label{vero1lcfsun}
f(\by | \mu, \sigma^{2}, \bDelta)
= \frac{2^{Lm} \exp\left\{-\sum_{i=1}^{L}\frac{(\ln y_{i}-\mu)^{2}}{2\sigma^{2}}\right\} }
{(2\pi\sigma^{2})^{L/2}\prod_{i=1}^{L}y_{i} }
 \prod_{i=1}^{L} \Phi_{m}\left(\frac{\bDelta'(\ln y_{i}-\mu)}{\sigma}|\bDelta^*\right).
\end{equation}

As discussed in \cite{QuLoSi15}, if the population has a LCFUSN distribution,  it is not
easy to elicit a  prior distribution for the skewness parameter $\bDelta$ when it has a very
general structure. Let us  consider a more parsimonious model where $\bDelta = \delta {\bf{1}}_{m}^T$.
In this case, the matrix ${\bf{I}}_m-\bDelta' \bDelta$ is positive definite if  $\delta$
belongs to the interval $(-1/\sqrt{m}, 1/\sqrt{m})$. Also,
consider that {\it{a priori}} $\mu$, $\sigma$ and $\delta$ are independent and such that $\mu \sim N(\mu_0,v)$,
$\sigma^{2} \sim IG(\alpha, \beta)$ and $\delta \sim U(-1/\sqrt{m}, 1/\sqrt{m})$, where $\mu_0 \in \RR$,  $v$, $\alpha$ and $\beta$
are non-negative numbers.  

The posterior distributions are easier obtained if we 
first apply the logarithmic transformation to data, that is,
if the original data is such that $Y_{i} \sim LCFUSN_{1,m}(\mu,\sigma^2,\delta{\bf{1}}_{m}^T)$,
then the transformed data is $Z_{i}=\ln Y_{i} \sim CFUSN_{1,m}(\mu,\sigma^2,\delta{\bf{1}}_{m}^T)$.
After that, consider the  stochastic representation of the CFUSN family (see \cite{ArGe05}) in terms of convolutions, which establishes that
$Z_{i}\buildrel d \over = \gamma {\bf{1}}_{m}^T|\bX_{i}|+[\tau^2]^{1/2}V_{i} +\mu,
$
where $\gamma=\delta \sigma$, $\tau^2 = \sigma^2 - \gamma^2$, $\bX_{i} \sim N_{m}({\bf{0}}, {\bf I}_{m})$, $V_{i} \sim N({0}, 1)$,
$\bX_{i}$ and $V_{i}$
are independent random quantities and $|\bX_{i}|=(|X_{i1}|,...,|X_{im}|)'$. 
Now, we can  hierarchically
represent our model as
\allowdisplaybreaks
\begin{eqnarray}\label{hier}
 Y_{i}  = \exp (Z_{i}); \,\,\,\, Z_{i}|\bX_{i}=\bx_{i}\sim N( \mu+ \gamma{\bf{1}}_{m}^T |\bx_{i}|,
 \tau^2);\,\,\,\, \bX_{i} \sim N_{m}({\bf{0}}, {\bI}_{m}).
\end{eqnarray}
Based on this stochastic representation, the variables $\bX_i$ can be considered latent variables in our model and the estimates
of $\mu,$ $\sigma^2,$ $\delta$ and $\bX_1,\dots,\bX_n$ can be obtained using the following likelihood function 
\allowdisplaybreaks
$$
f({\bf{z}}\mid {\bf{x}},\mu,\sigma^2,\delta)=\prod_{i=1}^{n}\phi(z_i\mid\mu+ \gamma{\bf{1}}_{m}^T |\bx_{i}|,
 \tau^2).
$$
Under this model representation the full conditional distributions (fcd) for the parameter $\mu$, $\sigma^2$ and $\delta$
and for the latent vector $\bX_i$, $i=1, \dots, L$ are, respectively,
\allowdisplaybreaks
\begin{eqnarray*}
\mu \mid  \sigma^2, \delta, \bZ, \bX &\sim& N\left( \frac{\mu_0\tau^2 + v\sum_{i=1}^{L}(z_i-\gamma \sum_{j=1}^{m}|x_{ij}|)}{Lv +\tau^2}, \frac{\tau^2v}{Lv +\tau^2}\right),\;\;\;\;\;\;\; \\[2mm]
f(\sigma^2 \mid  \mu, \delta, \bZ, \bX) &\propto&  \left( \frac{1}{\sigma^2}\right)^{\alpha + L/2+1}
    \exp\left\{ \frac{\sum_{i=1}^{L} (z_i - \mu) \delta \sum_{j=1}^{m}|x_{ij}|}{ \sigma(1-\delta^2)} \right\}  \\
    &\times& \exp\left\{- \frac{2 \beta (1-\delta^2) + \sum_{i=1}^{L} 2(z_i - \mu)^2}{2 \sigma^2(1-\delta^2)}\right\},\\[2mm]
f(\delta \mid \mu, \sigma^2, \bZ, \bX) &\propto&   \exp\left\{\frac{-\sum_{i=1}^{L}(z_i -\mu -\sigma\delta \sum_{j=1}^{m} |x_{ij}| )^2}{2\sigma^2(1 - \delta^2)} \right\}\\[2mm]
&\times&  \frac{1}{[1-\delta^2]^{L/2}} {\bf{1}}\{\delta \in (-1/\sqrt{m}, 1/\sqrt{m})\}, \\[2mm]
f(\bX_i \mid \mu, \sigma^2, \delta, \bZ, \bX_{(-i)} ) &\propto&
\exp\left\{ \sum_{i=1}^{L}\left[(z_i - \mu -\sigma\delta \sum_{j=1}^{m}|x_{ij}|)^2 - \frac{\sum_{j=1}^{m}x_{ij}^2}{2}\right] \right\}.
\end{eqnarray*}

The Gibbs sampler  can be used to sample from the posterior fcd of $\mu$.
The posterior fcd of $\sigma$, $\delta$ and $ \bX_i$, $i=1, \dots, n$, do not have
closed forms and the Metropolis-Hastings algorithm can be used to sample from such distributions.
Alternatively, we can assume that  $\mu$, $\gamma$ and $\tau^2$ are independent with  $\mu \sim N(\mu_0,v)$,
$\gamma \sim N(m, W)$ and $\tau^2 \sim IG(a,d)$. By assuming this,  the fcd of $\mu$ and $\bX_i$ remains the same
as before and of $\gamma$ and $\tau^2$ are, respectively,
\begin{eqnarray*}
\tau^2 \mid \mu, \gamma,  \bZ, \bX &\sim& IG \left( a+ \sum_{i=1}^{L}(z_i - \mu- \gamma{\bf{1}}_{m}^T\mid \bx_i\mid)^2, d+L\right)\\
\gamma\mid \mu,\tau^2,  \bZ, \bX &\sim&  N \left( \frac{m\tau^2 + W \sum_{i=1}^{L}(z_i + \mu){\bf{1}}_{m}^T\mid \bx_i\mid}{W^*},\frac{\tau^2W}{W^*} \right),
\end{eqnarray*}
where $W^*= \tau^2 +W \sum_{i=1}^{L}{\bf{1}}_{m}^T\mid \bx_i\mid$. This strategy facilitates the
implementation of the MCMC and may help its convergence. However, we lose the interpretation of the parameters
which can make the elicitation of prior distributions a more hard task.
Moreover, the hierarchical representation in Equation
$(\ref{hier})$ allows us to use  WinBUGS to obtain samples from the posteriori distributions.
For a more detailed discussion on Bayesian inference in the LCFUSN family see \cite{QuLoSi15}.

For each sample $(\mu^t, (\sigma^2)^t, \delta^t)$   of the posterior distribution of $(\mu, \sigma^2, \delta)$
we can obtain a sample of the posterior of the entropy $H_Y(\mu, \sigma^2, \delta)$. Posterior summaries of  $H_Y(\mu, \sigma^2, \delta)$
such as means, modes and HPDs can be approximated in the usual way. Similar procedure can be used to
sample from the posterior distributions of the MI and the KL divergence.
Another way to estimate  $H_Y(\mu, \sigma^2, \delta)$  is to plug
the posterior point estimates (usually, posterior means or modes) in its expression.
One disadvantage of this procedure  is that the posterior uncertainty about the parameters
is not considered in the estimation of $H_Y(\mu, \sigma^2, \delta)$.

\section{Applications}
\label{Apl}

\subsection{Simulated data sets analysis}\label{simul}

We run a simulation study comparing the LCFUSN distribution  with the well established  LSN and LN distributions.
A sample of size 3000 is generated from each one of the following distributions: the $LCFUSN_{1,3}(2, 0.8, \delta \bf{1}_{3,1})$  assuming $\delta=0.577$ (Data 1) and $\delta=-0.500$ (Data 2), the $LSN(2, 0.8, 0.990)$  (Data 3) and $LN(2,0.8)$ (Data 4).

To analyze the data, we  fit four models  assuming that 
        $Y_i \mid \mu, \sigma^2, \bDelta \sim LCFUSN_{1,m}(\mu, \sigma^2, \delta{\bf{1}}_{m,1})$, $m=2,3$, 
        $Y_i \mid \mu, \sigma^2, \bDelta \sim LSN(\mu, \sigma^2, \delta)$ and 
        $Y_i \mid \mu, \sigma^2, \bDelta \sim LN(\mu, \sigma^2)$.
To complete the model specification, in all cases we  assume  flat prior distributions for all parameters by eliciting   $\mu \sim N(0,100)$,
$\sigma^2 \sim IG(0.1, 0.1)$ and $\delta \sim U(\frac{-1}{\sqrt m},\frac{1}{\sqrt m})$.

Table 1 shows the posterior means and variances  provided by all fitted models. Comparing the posterior estimates
with the true values,  the  $LCFUSN_{1,3}$ distribution provides the less biased estimates for Data 1 and Data 2. For Data 3
this is attained if the LSN distribution is fitted. In general, the posterior means provided by all four models are comparable  in all cases. 
For data coming from a distributions with extreme values for the shape parameter,
that is Data 1 and Data 3, the LN  significantly underestimates the variance.  We also noticed that
the variance is overestimated by  LSN in  Data 1, and by all four models in Data 4. That can be a problem if, for instance, 
the main interest lies on estimate the quantiles of the distribution. 
Table 1 also  shows DIC, SlnCPO and entropy for the four models. The true model is correctly selected by
the entropy in Data 1 and Data 2, by the SlnCPO in Data 2 and Data 3 and by DIC  in Data 2, disclosing
that the entropy can be useful for model comparison. 

\begin{table}[]
    \centering
    \caption{Estimates of mean and variance and model selection statistics.}
    \label{None}
    \footnotesize 
    \begin{tabular}{ccccccc}
        \hline   
        \multicolumn{1}{c}{} &  \multicolumn{2}{c}{Estimates} & \multicolumn{4}{c}{Model Selection} \\ 
        \cline{2-3}\cline{5-7}
           Model & Mean & Variance & & DIC & SlnCPO & Entropy \\ 
        \hline
         \multicolumn{1}{c}{} &  \multicolumn{5}{c}{Data 1}  \\ 
          \multicolumn{1}{c}{} &  \multicolumn{2}{c}{$\mbox{True mean}=25.343$}&\multicolumn{2}{c}{$\mbox{True variance}=230.408$}&\multicolumn{2}{c}{$\mbox{True entropy}=5.127$}  \\ 
        \cline{2-7}
        
           LSN  & $26.046$ & $315.958$ && $-16410$ & $-0.679$ & $3.765$ \\
          $LCFUSN_{1,2}$   & $25.697$ & $267.102$ && $-63590$ & $-0.677$ & $4.487$ \\
          $LCFUSN_{1,3}$  & $25.167$ & $224.821$ && $-51270$ & $-0.681$ & $5.139$ \\
          $LN$   & $25.624$ & $186.508$ && $1457$ & $-0.729$ & $3.852$ \\
        \hline
         \multicolumn{1}{c}{} &  \multicolumn{5}{c}{Data 2}  \\ 
          \multicolumn{1}{c}{} &  \multicolumn{2}{c}{$\mbox{True mean}=3.323$}&\multicolumn{2}{c}{$\mbox{True variance}=3.769$}&\multicolumn{2}{c}{$\mbox{True entropy}=3.315$}  \\ 
        \cline{2-7}
           
            LSN  & $3.345$ & $3.929$ && $-1031$ & $-0.924$ & $1.922$ \\
          $LCFUSN_{1,2}$   & $3.367$ & $3.956$ && $-9540$ & $-0.870$ & $2.610$ \\
          $LCFUSN_{1,3}$   & $3.332$ & $3.861$ && $-19690$ & $-0.870$ & $3.331$ \\
          $LN$   & $3.303$ & $4.336$ && $5244$ & $-0.874$ & $1.910$ \\
        \hline
         \multicolumn{1}{c}{} &  \multicolumn{5}{c}{Data 3}  \\ 
          \multicolumn{1}{c}{} &  \multicolumn{2}{c}{$\mbox{True mean}=15.919$}&\multicolumn{2}{c}{$\mbox{True variance}=111.464$}&\multicolumn{2}{c}{$\mbox{True entropy}=3.248$}  \\ 
        \cline{2-7}

            LSN  & $15.91$ & $112.71$ && $-51690$ & $-0.597$ & $3.248$ \\
          $LCFUSN_{1,2}$   & $15.55$ & $83.16$ && $-357100$ & $-0.609$ & $3.909$ \\
          $LCFUSN_{1,3}$   & $15.44$ & $70.79$ && $-521900$ & $-0.629$ & $4.605$ \\
          $LN$   & $15.50$ & $64.59$ & &$1406$ & $-0.702$ & $3.326$ \\
         \hline
         
          \multicolumn{1}{c}{} &  \multicolumn{5}{c}{Data 4}  \\ 
          \multicolumn{1}{c}{} &  \multicolumn{2}{c}{$\mbox{True mean}=10.167$}&\multicolumn{2}{c}{$\mbox{True variance}=93.417$}&\multicolumn{2}{c}{$\mbox{True entropy}=3.187$}  \\ 
        \cline{2-7}

            LSN  & $10.43$ & $104.39$ && $6900$ & $-1.200$ & $3.235$ \\
          $LCFUSN_{1,2}$   & $10.74$ & $118.49$ && $6745$ & $-1.200$ & $3.955$ \\
          $LCFUSN_{1,3}$   & $10.73$ & $125.97$ && $6441$ & $-1.200$ & $4.664$ \\
          $LN$   & $10.29$ & $106.42$ & &$7201$ & $-1.200$ & $3.206$ \\
         \hline

    \end{tabular}
\end{table}

\begin{figure}[h]
  \centering
  \includegraphics[width=5.5cm]{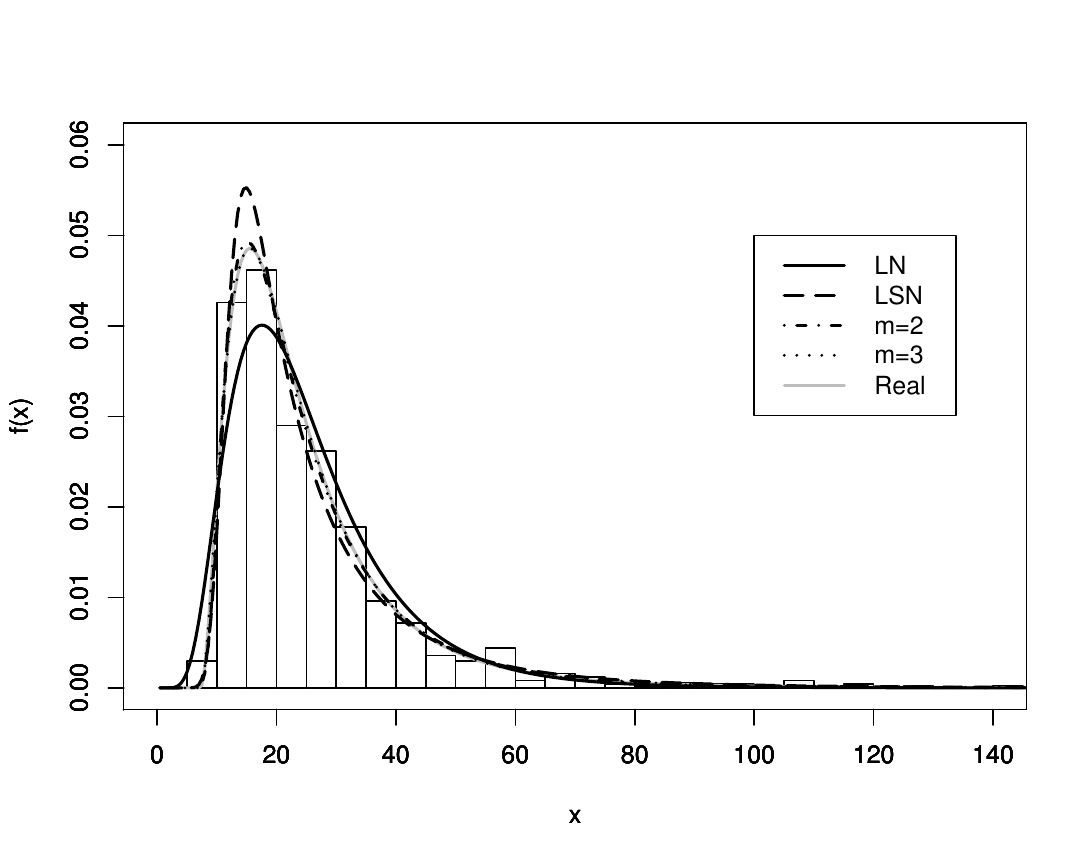}
  \includegraphics[width=5.15cm]{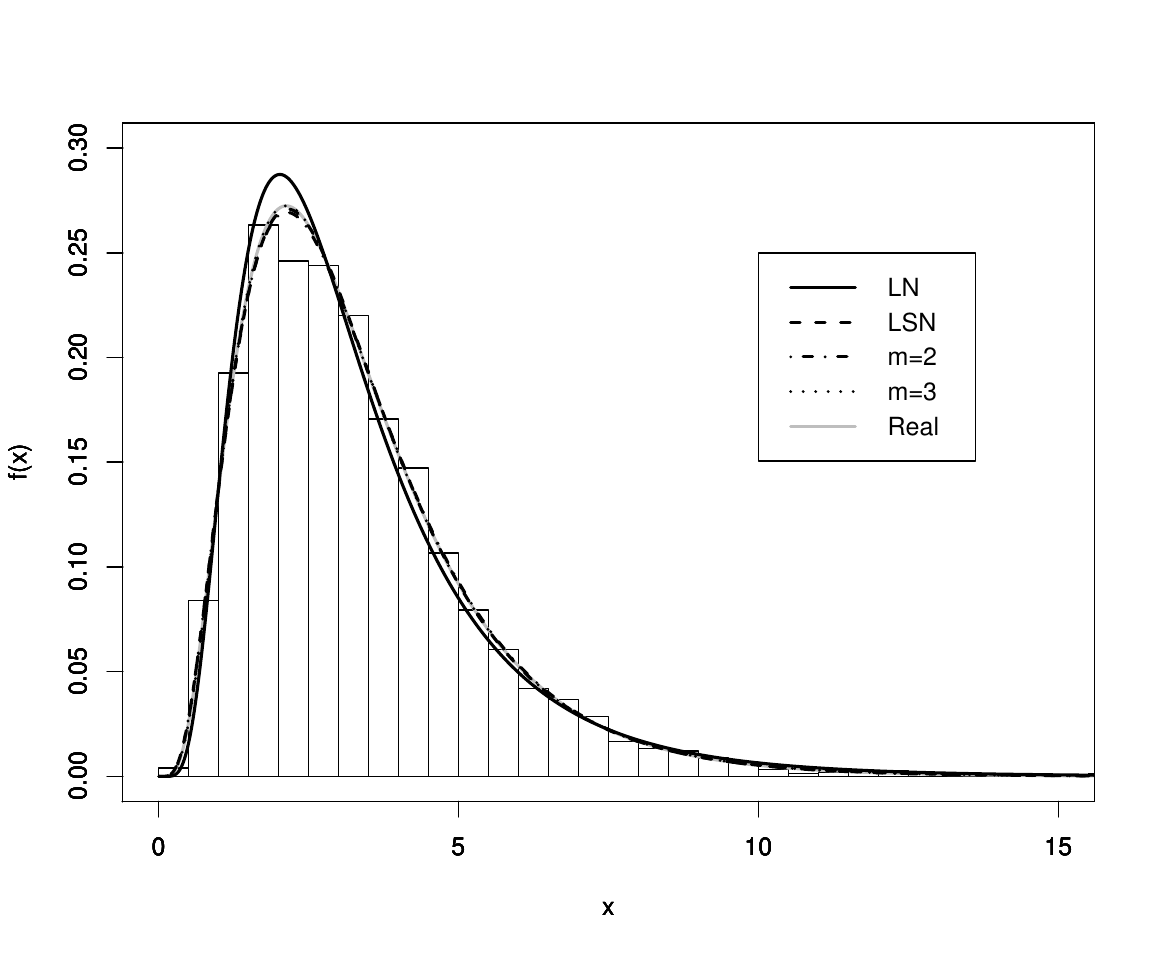}\\
  \includegraphics[width=5.15cm]{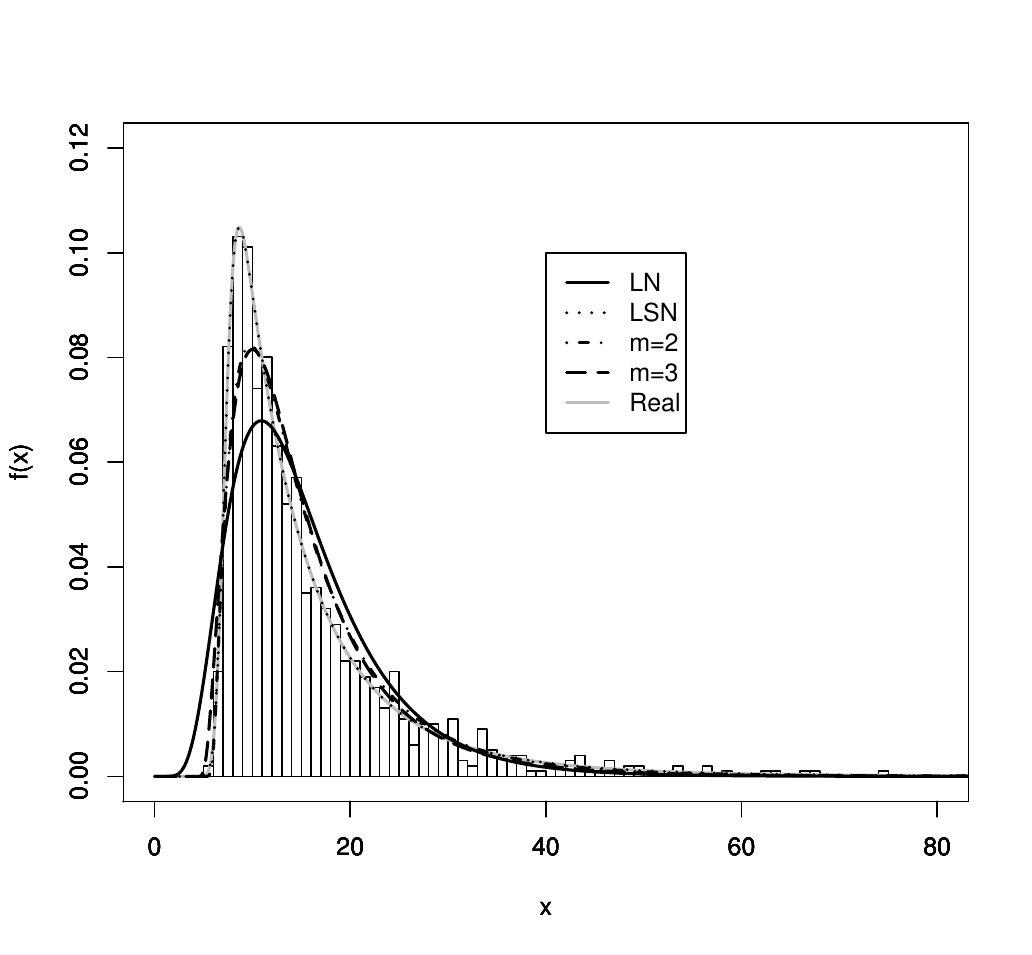}
  \includegraphics[width=5.5cm]{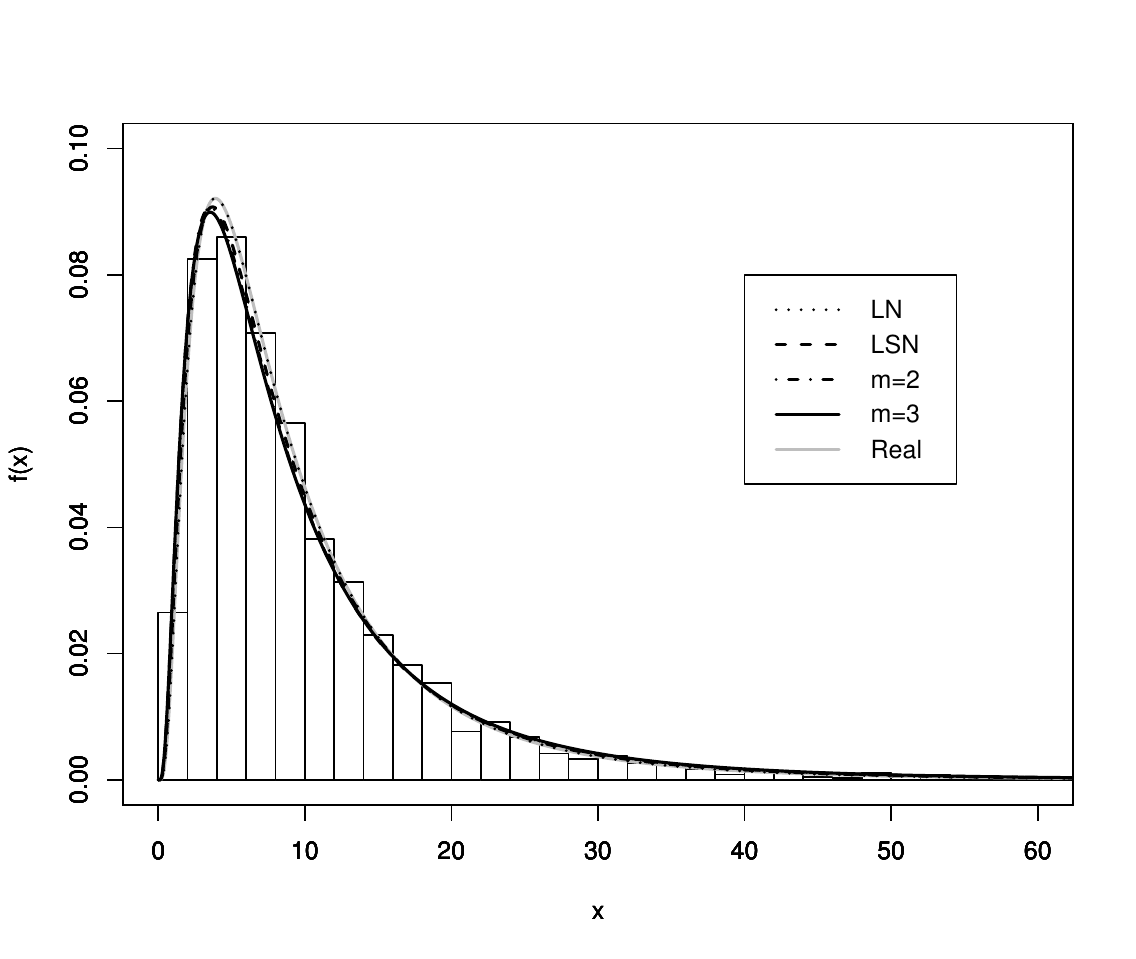}
  \caption{Estimated densities of $LN$, $LSN$, $LCFUSN_{1,2}$ and $LCFUSN_{1,3}$ models for Data 1(top left), Data 2 (top right), Data 3 (bottom left) and Data 4 (bottom right).}
\label{FiDist}  
\end{figure}

Figure \ref{FiDist} shows that, in Data 1 and Data 3,  the fitted
LN distribution  poorly estimates the left tail of the distribution.
The height of the mode is not well estimated by the LN distribution in Data 1, Data 2 and Data 3 and,
by the LSN distribution in Data 1. The proposed models do not estimate well the height of the mode  in Data 3
but they do it much better than the LN distribution. LSN and the proposed model 
are comparable in Data 2 and Data 4 and they are comparable to the LN  in Data 4, showing the flexibility of the 
LCFUSN distribution. 

\subsection{Case Study 1: Selecting model using entropy}

Quoting \cite{Ja57}, pp. 623, {\it{`` $\dots$ in making inference on the basis of partial information we
must use that probability distribution which has the maximum entropy subject to whatever is known.
This is the only unbiased assignment we can make; to use any other would amount to arbitrary assumption
of information which by hypothesis we do not have."}}

\begin{figure}[h!]
\centering
\includegraphics[width=4.5cm]{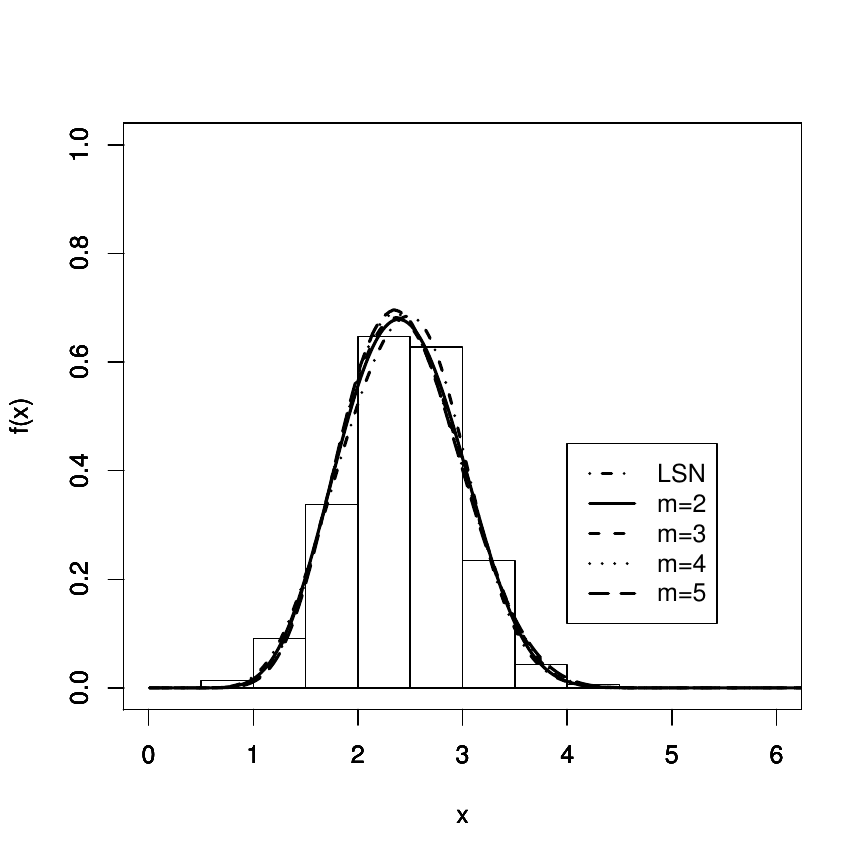}
\caption{Fitted LCFUSN and LSN  densities, precipitation data.}
\label{AplicPre}
\end{figure}

In this section we  analyze  the USA monthly precipitation data recorded from 1895 to 2007
by fitting LCFUSN distributions with different values for $m$. This data
is available at the National Climatic Data Center (NCDC) and consists of 1,344 observations
of the US precipitation index (PCL). Our main goal here is to consider  the Jaynes's
principle to select the best model, that is, we select the model that maximizes the  entropy. Models are also chosen
using two well-known tools for model selection, the conditional predictive ordinate (SlnCPO)
and  the deviance information criterion (DIC).
Denote by $Y_i$ the precipitation index in the $i$th month.
In \cite{QuLoSi15} the authors fitted the LSN (say, the LCFUSN with $m=1$)
and different  LCFUSN distributions and evaluate the gain
in assuming a higher dimensional skewing function to analyze the data.
We assume the same models, that is, we consider that
$Y_i \mid \mu, \sigma^2, \delta \sim LCFUSN_{1,m}(\mu, \sigma^2, \delta{\bf{1}}_{m,1})$
and postulate for all parameters the same flat prior distributions elicited in Subsection \ref{simul}.
We also let  $m$ to vary from $m=1$ (LSN) to $m=5$.
We name $M_i$ the model for which we assume $m=i$.
By considering such specifications and assuming the posterior means, in \cite{QuLoSi15}
it is obtained very close plug-in estimates of the true
density for all $m$  as can be noticed in Figure \ref{AplicPre}.

Table 2 shows the 95\% HPD, SlnCPO, DIC and the entropy
comparing all models. The entropy was
computed using the two different approaches presented in Section \ref{SeBaE}.
\begin{table}[h]
\centering
\label{roger10}
\caption{Model selection statistics, Precipitation data.}

\begin{tabular}{c|c|c|c|c|c|c}\hline


Model      & $E(\delta\mid \bY)$ & Mean Entropy & $95\%$ HPD & Plug-in Entropy & DIC & SlnCPO \\ \hline
$LN$ & $-$ &  $0.883$        &$[0.844;0.924]$ &  $0.890$ & $62.69$ & $-0.8815$ \\  \hline
$LSN$ & $-0.947$ &  $1.143$        &$[0.963;1.313]$ &  $1.158$ & $-13.19$ & $-0.83766$ \\  \hline
$m=2$      & $-0.686$ &  $1.164$        &$[1.101;1.236]$ &  $1.230$ & $-36.96$ & $-0.83545$ \\  \hline
$m=3$      & $-0.570$ &  $1.256$        &$[1.160;1.348]$ &  $1.314$ & $-112.40$& $-0.83765$\\  \hline
$m=4$      & $-0.497$ &  $1.467$        &$[1.364;1.566]$ &  $1.539$ & $-321.10$& $-0.84144$\\  \hline
$m=5$      & $-0.446$ & $1.777$         &$[1.675;1.880]$ &  $1.815$ & $-895.30$& $-0.81057$\\  \hline
\end{tabular}
\end{table}

\noindent We notice that the mean entropy increases with the
complexity of the model. Moreover, the HPDs for models $LSN$ and  $M_2$
point out that the amount of information
in our system is quite similar and this same conclusion holds for models $M_2$ and $M_3$.
The HDPs for models $M_4$ and $M_5$ do not intersect,
indicating that the entropy of these two models are  significantly different. 
By using the maximum entropy principle we decide for $M_5$ as the best model, and this decision is the same we make using
standard procedures for model selection, such as
DIC and the SlnCPO.
Moreover, we notice that, at least in this example,
the  entropy and the DIC have similar behavior
and leads to the same decision. It should be observed that the $LN$ model, which is the most common choice for this type of data analysis, provides the poorest fit among all models, as shown by the model selection statistics DIC and SlnCPO as well as by the Shannon entropy.


\subsection{Case Study 2: Clustering with KL divergence}

\begin{figure}[h!]
\centering
\includegraphics[width=4.5cm]{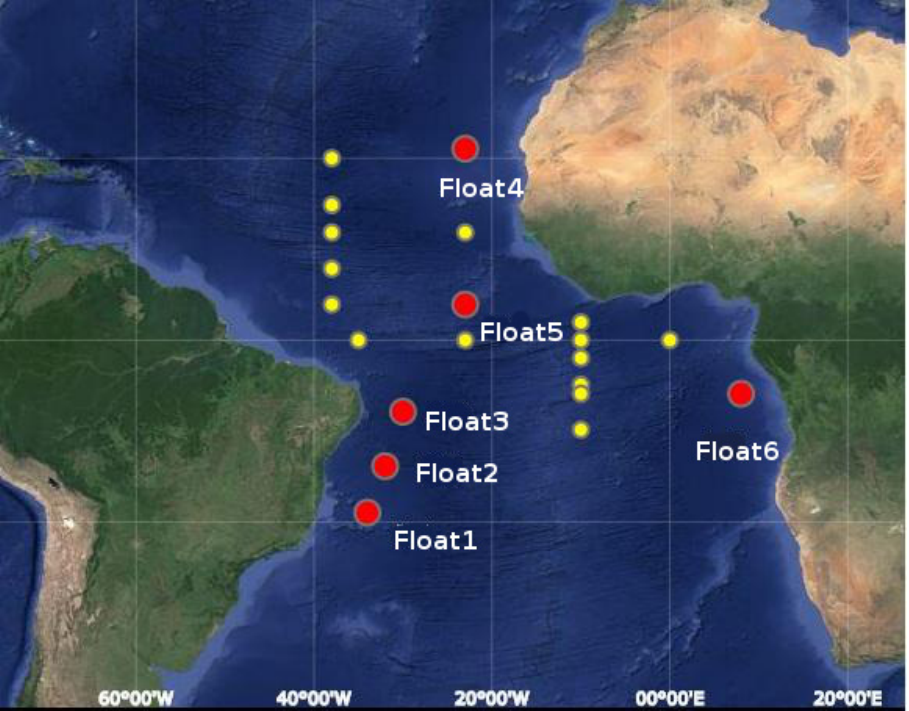}
\caption{Selected floats in Atlantic ocean.}
\label{FigBo}
\end{figure}

Climate in some Brazilian areas are  directly
influenced by some features of Atlantic ocean
such as surface temperature and humidity. To
better understand how such influence occurs,
a system of 21 floats  were installed over different
regions of Atlantic ocean  (see Figure \ref{FigBo})
and some of such features are daily measured.

We  choose six of this floats and consider the daily  relative humidity 
from September 11th, 1997 to September 22nd, 2014.  
Such floats (their locations are in parenthesis) are  named Float1 (19S 34W), Float2 (4S 32W), Float3(8S 30W), Float4(21N 23W),
Float5 (4N 23W) and Float6 (6S 8E) and are the big balls in Figure \ref{FigBo}.
Data are available at URL:http://goosbrasil.org/.
Our goal is to verify if the relative humidity level  in all these regions
have similar behavior. The importance of such kind of study from the statistical point of view
is to permit us to cluster information in similar floats looking
for improvements in the estimates. On the other hand, from the economic perspective,
it can be considered to define an optimal design for the system used to measure
information in Atlantic area excluding some of the floats which brings
similar information.

We assume two different distributions to model the humidity $Y$ in each float.
In the first case we consider that $Y_i \mid \mu, \sigma^2 \sim LN(\mu, \sigma^2)$,
which is a standard assumption to analyze this data.
As an alternative we also assume that $Y_i \mid \mu, \sigma^2, \bDelta \sim LCFUSN_{1,2}(\mu, \sigma^2, \delta{\bf{1}}_{1,2})$.
As in Case Study 1, flat prior distributions for the parameters are considered.

Table 3 provides the DIC comparing the fitted models
for each float as well as  the  KL divergence between the two distributions under consideration.
Considering the method introduced in \cite{Mc89} and assuming
that distributions are similar whenever the  KL divergence
is up to the cut point $0.14$,
we concluded 
that the LN and LCFUSN are very similar
to describe the humidity behavior in each float.
Despite of this,   the DIC points out that the
LCFUSN distribution is a better model for all floats.
Therefore, in the following, our analysis is based on
the LCFUSN distribution only.

\begin{table}[h]
\centering
	\label{TaMoSe2}
	\caption{DIC and Kullback-Leibler divergence between LCFUSN and LN distributions for each float.}

		\begin{tabular}{c|c|c|c|c}\hline
			Float & DIC-LN & DIC-LCFUSN & $D(LCFUSN||LN)$ & $D(LN||LCFUSN)$ \\
	
			1   &   $-6,310$ & $ -20,500$ & $0.0047$ & $0.0040$\\ \hline
			2   &   $-9,703$ & $-17,930$  & $0.0017$ & $0.0039$\\ \hline
			3   &   $-9,708$ & $-15,070$  & $0.0007$ & $0.0010$\\ \hline
			4   &   $-4,859$ & $ -104,100$& $0.0365$ & $0.0806$\\ \hline
			5   &   $-8,079$ & $-20,790$  & $0.0064$ & $0.0055$\\ \hline
			6   &   $-2,530$ & $-25,120$  & $0.0307$ & $0.0663$\\ \hline
		
		\end{tabular}
	\end{table}

Table 4 shows the KL divergence
comparing the LCFUSN distributions $f_i$ and $f_j$ of floats $i$ and $j$.
Since the KL divergence is asymmetric,
we will assume that data of different  floats have the same distribution,
and thus could be clustered, whenever $D(f_i|f_j)$  and $D(f_j|f_i)$
are both up to $0.14$.   It can be noticed from Table 4
that the relative humidity  measured on floats 3 and 6 are similar, hence one of them could be excluded from the system  without losing substantial information. In fact,  the predictive distributions based on data measured on Floats 3 and 6 
are quite similar to that obtained when we merge the data from these floats.
\begin{table}[h]
\centering
\label{TaClu}
\caption{Kullback-Leibler divergence between $f_i$ and $f_j$.}
\begin{tabular}{c|c|c|c|c|c|c} \hline

$Float$ & 1 & 2 & 3& 4 &5  &6\\ \hline

1   & $-$      & $0.2534$ & $0.5427$ & $0.1446$ & $1.0543$ & $0.6560$\\ \hline
2   & $0.1390$ & $-$      & $0.1646$ & $0.3346$ & $0.8198$ & $0.4549$\\ \hline
3   & $0.2433$ & $0.1602$ & $-$      & $0.2471$ & $0.3117$ & ${\bf{0.1249}}$\\ \hline
4   & $0.1761$ & $0.6623$ & $0.5055$ & $-$      & $0.4862$ & $0.3795$\\ \hline
5   & $0.6707$ & $1.0787$ & $0.4345$ & $0.2134$ & $-$      & $0.1766$\\ \hline
6   & $0.4235$ & $0.5416$ & ${\bf{0.1323}}$ & $0.2033$ & $0.0937$ & $-$\\ \hline

\end{tabular}
\end{table}

Table 5 shows the mean and the variance
of the prior predictive distribution for the relative humidity in all these cases.
As can be noticed, the predictive distributions have very close means and
similar variance.  However, when we merged the data, the variability on the humidity
on Float 6 is underestimated by a small amount. That probably happens
because data in Float 6 discloses a high degree of skewness as
can be noticed from Table 6.  It is also noteworthy that 
the merged  data  better capture the features of data measured
on Float 3.

\begin{table}[h]
\centering
\label{TaPre}
\caption{Prior predictive mean and variance based on data in Floats 3 and 6 and the merged.}
\begin{tabular}{c|c|c}\hline
 $$& $E(Y)$ & $Var(Y)$\\ \hline
Float 3 &  $78.66$ & $14.91$\\ \hline 
Float 6 &  $80.45$ & $17.20$\\ \hline
Merged  &  $79.56$ & $15.41$\\ \hline

\end{tabular}
\end{table}

\begin{table}[h]
\centering
\label{TaPoSu}
\caption{Posterior summaries data of  Floats 3 and 6 and the merged data.}
\begin{tabular}{cccccccc}\hline
\multicolumn{1}{c}{}& \multicolumn{2}{c}{$\mu$} &\multicolumn{2}{c}{$\sigma$} &\multicolumn{3}{c}{$\delta$} \\
\cline{2-3}\cline{4-5}  \cline{6-8}
$$ & Mean & St. Dev. &  Mean & St. Dev. &  Mean & St. Dev. & $95\%$HPD  \\ 
Float 3 &  $4.410$ & $0.005$ &  $0.059$ &  $0.002$ & $-0.493$ & $0.037$ & $[-0.564,-0.428]$\\ 
Float 6 &  $4.477$ & $0.004$ &  $0.083$ &  $0.003$ & $-0.687$ & $0.007$ & $[-0.699,-0.672]$\\
Merged  &  $4.420$ & $0.005$ &  $0.062$ &  $0.002$ & $-0.521$ & $0.031$ & $[-0.575,-0.457]$\\ \hline

\end{tabular}
\end{table}

\section{Final Comments}

In this paper we presented some results related to the 
entropy and KL divergence  of
two classes of multivariate distributions recently
introduced in the literature, the multivariate log-canonical
fundamental skew-normal (LCFUSN) and the canonical fundamental
skew-normal (CFUSN) distributions.  We also obtained
the MI for distributions in  these families.
Bayesian approach is considered to estimate the  entropy and
KL divergence. Such measures were computed in two
different ways, obtaining their posterior distribution and
using the plug-in method where the posterior means of
the parameters were assumed as estimators.
To illustrate the use of some results, entropy was used for model comparison.
We concluded that Jaynes's principle selected the
same model as DIC and CPO.
KL divergence was  used to compare the distribution
of relative humidity collected in scattered floats on different
regions of Atlantic ocean. We concluded 
that the humidity in some floats has a very similar behavior
and some of such floats could be removed from the
system if the only goal were to measure the humidity.

An interesting topic for future research is to consider some generalizations of Shannon entropy, such as 
Mathai's generalized entropy and R\'enyi's entropy, in the $LCFUSN$ and $CFUSN$ families of distributions.
Mathai's generalized entropy was used in \cite{JN08}, who proved that the pathway model can be obtained by optimizing such entropy.
\section*{Acknowledgements}
The authors would like to thank the Editor, the Associate Editor and the referees for their many helpful comments and suggestions. The authors also thank Professor Sacha Friedly for his suggestions. The research of Marina de Queiroz was partially supported by  CAPES (Coordena\c{c}\~ao de Aperfei\c{c}oamento de Pessoal de N\'ivel Superior) and CNPq (Conselho Nacional de Desenvolvimento Cient\'ifico e Tecnol\'ogico). Rosangela Loschi would like to thank to CNPq, [grant number 301393/2013-3], [grant number 306085/2009-7], for a partial allowance to her researches. The research of Roger Silva was partially supported by FAPEMIG [grant number APQ-02743-14].

\section{Appendix}

In this section we present the proofs of propositions that appear in the text.\\

{\textbf {Proof of Proposition \ref{enlsncap2}:}} First note that if $\bX\sim CFUSN(\bmu,\bf{\Sigma},\bf\Delta)$, then $\bX=\bmu + \bf{\Sigma}^{1/2}\bZ$, where $\bZ\sim CFUSN(\bf\Delta)$,
and $H_{\bX}=\frac{1}{2}\ln|\bf{\Sigma}|+H_{\bZ}$. Thus we only need to calculate $H_{\bZ}$. Consider the pdf given in Equation (\ref{CFUSN}).
Thus, we have
$H_{CFUSN(\bf{\Delta})}=-E(\ln \phi(Z)) - E_{\bZ}[\ln(2^{m}\Phi_{m}(\bDelta'\bZ|\bDelta^{*}))]=
\frac{n}{2} \ln 2\pi  + \frac{1}{2} \sum_{i=1}^{n}E(Z_{i}^{2})-E_{\bZ}[\ln(2^{m}\Phi_{m}(\bDelta'\bZ|\bDelta^{*}))].$
The proof is concluded by noticing that $\sum_{i=1}^{n}E(Z_{i}^{2})=E(\bZ'\bZ)=\bmu_0\bmu_0+tr(\bSigma_0)$, where
$\bmu_0=E(\bZ)$ and $\bSigma_0=Var(\bZ)$.\\

{\textbf{Proof of Proposition \ref{propentropiarelativa}(ii):}}
By definition it follows that
\begin{equation}\label{entrorelativa2}
 D(f_{\bY}||f_{\bZ}) = - H_{CFUSN_{n,m_{2}}(\bmu_{2}, \mathbf{\Sigma}_{2}, \mathbf{\Delta_{2}})}- \int_{\mathbb{R}^{n}} f_{\bY}(\by) \ln [f_{\bZ}(\by)]d\by.
\end{equation}
The integral at the right side of Equation (\ref{entrorelativa2}) can be calculated as follows:
\begin{eqnarray}
\label{inf3}
&&\int_{\mathbb{R}^{n}} f_{\bY}(\by) \ln [f_{\bZ}(\by)]d\by = -\sum_{i=1}^{n} E_{\bY}(Y_{i}){-\frac{1}{2}}\ln |\bSigma_{1}|-\frac{n}{2}\ln 2\pi \nonumber\\&-& \frac{1}{2}E_{\bY}\left((\bY -\bmu_{1})'\bSigma_{1}^{-1}(\bY -\bmu_{1})\right)
 + E_{\bY}\left(\ln\left[2^{m_{1}} \Phi_{m_{1}}(\bDelta_{1}'\bSigma_{1}^{-\frac{1}{2}}( \bY -\bmu_{1})|\bDelta_{1}^{*})\right]\right),
\end{eqnarray}
where $\bY \sim CFUSN _{n,m_{2}}(\bmu_{2},\mathbf{\Sigma}_{2}, \bDelta_{2})$ and $Y_{i}$ is the $i$-th component of $\bY$.
Replacing $(\ref{inf3})$ in $(\ref{entrorelativa2})$ we have that
\allowdisplaybreaks
\begin{eqnarray*}
&&D(f_{\bY}||f_{\bZ})= \frac{1}{2} \ln\left[\frac{\bSigma_{1}}{\bSigma_{2}} \right]+E_{\bY}\left( \ln \left[\frac{2^{m_{2}-m_{1}}\Phi_{m_{2}}(\bDelta_{2}'\bSigma_{2}^{-\frac{1}{2}}( \bY -\bmu_{2})|\bDelta_{2}^{*})}{\Phi_{m_{1}}(\bDelta_{1}'\bSigma_{1}^{-\frac{1}{2}}( \bY -\bmu_{1})|\bDelta_{1}^{*})} \right]\right)\nonumber \\
&&+ \frac{1}{2}E_{\bY}\left((\bY -\bmu_{1})'\bSigma_{1}^{-1}(\bY -\bmu_{1})\right)-\frac{E_{\bY_0}(\bY_{0}'\bY_0)}{2},
\end{eqnarray*} where $\bY_0\sim CFUSN_{n,m}(0,{\bf I}_n,\bDelta)$. The proof follows from Equation (\ref{EspCFSUN})
by noticing that $E_{\bY_0}(\bY_{0}'\bY_0)=\bmu_0\bmu_0+tr(Var(\bY_0))$, where $\bmu_0=E(\bY_0)$.


\begin{thebibliography}{15}
\bibitem{AhGo89}
Ahmed, N. A. \& Gokhale, D. V. (1989). Entropy expressions and their estimators for multivariate
distributions. {\emph{IEEE Trans. Inform. Theory}}, 35, 688-692.


\bibitem{ArCoGe11}
Arellano-Valle, R. B., Contreras-Reyes, J. E. \& Genton, M. G. (2013). {Shannon entropy and mutual information
for multivariate skew-elliptical distributions.} {\em Scandinavian Journal of Statistics},
40, 42-62.

\bibitem{ArGe05}
Arellano-Valle, R. B. \& Genton, M. G. (2005). On fundamental skew distributions.
{\em Journal of Multivariate Analysis},  96(1), 93-116.


\bibitem{Az85}
Azzalini, A. (1985). A class of distributions which includes the normal ones.
{\em Scandinavian Journal of Statistics}, 12, 171-178.


\bibitem{AzCaKo03}
Azzalini, A., Dal Cappello, T. \& Kotz, S. (2002). {Log-skew-normal and log-skew-t distributions as models for family income data.}
{\em Journal of Income Distribution}, 11(3,4), 12-20.

\bibitem{AzDa96}
Azzalini, A. \& Dalla Valle, A. (1996). { The multivariate skew-normal distribution.}
{\em Biometrika}, 83, 715-726.

\bibitem{Be79}
Bernardo, J. M. (1979). Reference  posterior distributions for Bayesian inference (with discussion).
{\em Journal of the Royal Statistical Society, Series B}, 41, 113-147.


\bibitem{CoAr12}
Contreras-Reyes, J. E. \& Arellano-Valle, R. B. (2012). {Kullback-Leibler divergence measure for
multivariate skew-normal distributions.} {\em Entropy}, 14(9), 1606-1626.

\bibitem{CoTh06}
Cover, T. M. and  Thomas, J. A. (2006). {\em Elements of information theory,} 2nd ed. Wiley, New Jersey.




\bibitem{GuSr10}
Gupta, M. \& Srivastava, S. (2010). {Parametric Bayesian estimation of differential
entropy and relative entropy} {\em Entropy}, 12, 818--843.

\bibitem{GuMu04}
Guti\'errez-Pe\~na, E. \& Muliere, P. (2004). {Conjugate priors represent strong pre-experimental assumptions.}
{\em Scandinavian Journal of Statistics}, 31, 2325-246.

\bibitem{JaGu08}
Javier, W. R. \& Gupta, A. K. (2008). {Mutual information for the mixture of two multivariate normal
distributions}. {\em Far East Journal of Theoretical Statistics},  26, 47-58.

\bibitem{JaGu09}
Javier, W. R. \& Gupta, A. K. (2009). {Mutual information for certain multivariate distributions.}
{\em Far East Journal of Theoretical Statistics}, 29, 39-51.

\bibitem{Ja57}
Jaynes, E. T.  (1957). {Information Theory and Statistical Mechanics}.
{\em The Physical Review}, 106(4), 620-630.


\bibitem{Ja68}
Jaynes, E. T.  (1968). {Prior probabilities.}
{\em IEEE Transactions Systems, Science and Cybernetics}, 4, 227-291.

\bibitem{JN08}
Jose, K. K. \& Naik, S.R.  (2008). {A class of asymmetric pathway distributions and an entropy interpretation.}
{\em Physica A}, 387, 6943-6951.

\bibitem{Ku78}
Kullback, S. (1978). {\em Information theory and statistics,} Dover Edition, Gloucester.

\bibitem{KuLe51}
Kullback, S. \& Leibler, R. A. (1951). On information  and sufficiency. {\em  Annals of Mathematical Statistics}, 22, 79-86.

\bibitem{LaPeJoEn00}
Larsen K, Petersen JH, Budtz-J$\o$rgensen E and Endahl L. 
Interpreting parameters in the logistic regression model with random
  effects. {\em{Biometrics}}. 2000;  {\bf{56}}: 909--914.

\bibitem{MaGe10}
Marchenko, Y. V. \& Genton, M. G. (2010). Multivariate log-skew-elliptical distributions with
applications to precipitation data. {\em Environmetrics}, 21, 318-340.

\bibitem{Mc89}
McCulloch, R. E. (1989). Local model influence. {\em Journal of the American Statistical Associ-
ation}, 84(406), 473-478.

\bibitem{QuLoSi15}
Queiroz, M. M., Loschi, R. H. \& Silva, R. W. C. (2016). Multivariate Log-Skewed Distributions with normal
kernel and its Applications. {\em Statistics}, {\bf 50}(1), 157-175.

\bibitem{SaLoAr13}
Santos, C.C., Loschi, R.H. and Arellano-Valle, R.B. (2013).
Parameter Interpretation  in  Skewed  Logistic Regression
with Random Intercept. {\em Bayesian Analysis}, {\bf 8}(2), 381--410.

\bibitem{Sa02}
Sahu, S. K. (2002).
Bayesian Estimation and model selection choice in item response models.
{\em Journal of Statistical Computation and Simulation}, 72(3), 217-232.

\bibitem{Sh48}
Shannon, C. E. (1948).
A mathematical theory of communication. {\em The Bell System Technical  Journal}, 27, 379-423, 623-656.

\bibitem{YuCl99}
Yuan, A. \& Clarke, B. S. (1999). A minimally informative likelihood for decision analysis:
Illustration and robstness. {\em Canadian Journal of Statistics},  27, 649--665.
\end{thebibliography}
\end{document}